\documentclass [11pt]{amsart} 
\usepackage {amsmath,amssymb, amscd, amsthm, epsfig, color}
\input xy
\xyoption {all}

\newtheorem{theorem}{Theorem} 
\newtheorem {lemma}{Lemma} 
\newtheorem {corollary}{Corollary} 
\newtheorem {proposition}{Proposition}

\theoremstyle{definition}
\newtheorem{remark}{Remark}

\theoremstyle {definition}

\theoremstyle {definition}

\def \pr {{\mathbb {P}^r}}

\def \b {{\bullet}}
\def \one {{\{1\}}}
\def \two {{\{2\}}}
\def \fscheme {{\overline {M}_{0,n} (X, \beta)}}

\def \schemea {{\overline M_{0,1}(\pr, d)}}
\def \schemeb {{\overline M_{0,2}(\pr, d)}}
\def \gscheme {{\overline M_{0,n}(\bg, d)}}

\def \p {{\mathbb {P}}}
\def \o {{{\mathcal O}_{\mathbb P}}}
\def \mon {{M_{0,n}}}
\def \monbar {{\overline {M}_{0,n}}}

\def \bg {{{\bf G}}}
\def \s {{\mathcal {S}}}
\def \q {{\mathcal {Q}}}

\def \o {{\mathcal O_{\mathbb P^1}}}
\def \c {{\mathbb C}}
\def \cs {{\mathbb C^{\star}}}

\def \mpnpmbar {{\overline {M}_{0,0}(\mathbb P^n \times \mathbb P^m, (d,e))}}
\def \mpnpm {{M_{0,0}(\mathbb P^n \times \mathbb P^m, (d,e))}}

\begin {document}

\title[Divisors on the Moduli Spaces of Stable Maps]{Divisors on the Moduli Spaces of Stable Maps to Flag Varieties and Reconstruction}
\author {Dragos Oprea}
\address {Department of Mathematics,}
\address {Massachusetts Institute of Technology,}
\address {77 Massachusetts Avenue, Cambridge, MA 02139.}
\email {oprea@math.mit.edu}
\date{} 

\begin {abstract}

We determine generators for the codimension $1$ Chow group of the moduli spaces of genus zero stable maps to flag varieties $G/P$. In the case of $SL$ flags, we find all relations between our generators, showing that they essentially come from $\overline M_{0,n}$. In addition, we analyze the codimension $2$ classes on the moduli spaces of stable maps to Grassmannians and prove a new codimension $2$ relation. This will lead to a partial reconstruction theorem for the Grassmannian of $2$ planes. 

\end {abstract}
\vskip1in 

\maketitle
 
In \cite {O1} and \cite {O2} we began the study of the tautological rings of the moduli spaces of stable maps to general flag varieties. We conjectured that the cohomology of these moduli spaces is entirely tautological and proved this conjecture for $SL$ flags. In this note, we bring more evidence in favor of our conjecture, establishing it completely for general flag varieties in codimension one, using a different method. 

More precisely, let $X$ be a projective homogeneous space. The (coarse) moduli spaces $\fscheme$ parametrize marked stable maps to $X$ in the cohomology class $\beta \in H^{2}(X, \mathbb Z)$. These moduli spaces are related by a complicated system of natural morphisms which we enumerate below: \begin {itemize} 
\item \text {forgetful morphisms}: $\pi:\overline M_{0,S} (X, \beta) \to {\overline M}_{0, T}(X, \beta)$ \text {defined whenever} $T\subset S$. \item \text {gluing morphisms which produce maps with nodal domain curves}, $$gl:{\overline M}_{0, S_1\cup \{\star\}}(X, \beta_1) {\times_{X}} {\overline M}_{0, \{\bullet\} \cup S_2} (X, \beta_2) \rightarrow \overline M_{0, S_1\cup S_2}(X, \beta_1+\beta_2).$$ \item {\text evaluation morphisms to the target space}, $ev_i:\fscheme \to X$ for all $1\leq i\leq n$.\end {itemize}

The system of tautological classes on $\fscheme$ is defined as the smallest subring of $H^{\star}(\fscheme)$ (or of the Chow ring $A^{\star}(\fscheme)$) with the following properties: 
\begin {itemize}
\item The system is closed under pushforwards and pullbacks by the natural morphisms. 
\item All monomials in the evaluation classes $ev_i^{\star} \alpha$ for $\alpha \in H^{\star} (X)$ are in the system. 
\end {itemize}

Typical examples of tautological classes are the following. Let $\alpha_1, \ldots, \alpha_p$ be cohomology classes on $X$. The class $\kappa(\alpha_1, \ldots, \alpha_p)$ is defined as the pushforward via the forgetful projection $\pi:{\overline M}_{0, n+p}(X, \beta)\to \fscheme$: $$\kappa(\alpha_1, \ldots, \alpha_p)=\pi_{\star}(ev_{n+1}^{\star} \alpha_1 \cdot \ldots \cdot ev_{n+p}^{\star} \alpha_p).$$ 

One of the main results of this note is the following: 

\begin {proposition} \label {gmodpa} Let $X$ be any projective homogeneous space. All (rational coefficient) complex codimension $1$ classes on $\fscheme$ are tautological. An explicit system of generators can be written down. \end {proposition}

We obtain sharper results for $SL$ flag varieties. These results generalize Pandharipande's (\cite {divisors}) who worked out the case of maps to projective spaces. However, Pandharipande's argument does not extend to this setting, essentially because there is no convenient description of the space of morphisms to flag varieties which is suitable for computations. Instead, we propose a different approach: the previous proposition gives generators for the codimension $1$ rational Chow group. We cut down the number of generators by exhibiting relations between the various $\kappa$ classes on $\fscheme$. We refer the reader to lemma $\ref {relation}$ for the precise equations. We show independence of the remaining classes by a dimension calculation. This calculation involves two ingredients: localization and a computation of the symmetric group invariants on the second cohomology of moduli spaces of marked rational curves. This latter computation is similar to Getzler's (\cite {zero}); the exact result is contained in lemma $\ref {dim}$. The argument will enumerate all the $\cs$ fixed loci on $\fscheme$ which have at most one negative weight on their normal bundle. 

To explain our result in more detail, we let $X$ be the $SL$ flag variety parameterizing $l$ successive quotients of fixed dimensions of a vector space $V$. We consider the tautological sequence on $X$: $$V\otimes \mathcal O_{X} \to \q_1\to\ldots \to \q_l\to 0.$$ We will prove that:

\begin {theorem}\label{diva} The following classes span the second cohomology group/codimension $1$ Chow group with rational coefficients of $\fscheme$: \begin {itemize} 
\item the boundary divisors,
\item the classes $\kappa(c_1(\q_i)^2)$ for all $1\leq i\leq l$ and the {\bf nonzero} classes $\kappa(c_2(\mathcal K_i))$ for $0\leq i\leq l$. Here $\mathcal K_i$ is the kernel of $\q_i \to \q_{i+1}$, and, by convention, $\q_0=V\otimes \mathcal O_{X}$ and $\q_{l+1}=0$.
\item when $n=1$ or $n=2$, we add any one of the evaluation classes $ev_i^{\star} c_1(\q_j)$.
\end {itemize}
The class $$\sum_{i} \kappa(c_2(\mathcal K_i))+ \sum_{i} \left(\frac{d_{i-1}+d_{i+1}}{2d_i}-1\right) \kappa(c_1(\mathcal Q_i)^2)$$ is sum of boundary classes. All other relations between these generators come from $\overline M_{0,n}$. The dimension of $H^2(\fscheme)$ is given by $$\left [2^{n-1}(d_1+1) \ldots (d_l+1)+\frac{1}{2}\right] - 1 - \binom {n}{2}+h^{4}(X) - \binom {h^2(X)}{2} .$$ Here $h^2(X), h^4(X)$ are the second and fourth Betti numbers of $X$,  $d_i=c_1(\q_i)\cdot \beta$, and $\left[\;\right]$ denotes the integer part.
\end {theorem}

Pandharipande's arguments (\cite {divisors}) show that all top intersection numbers of these divisors can be expressed as certain $n$ point genus $0$ primary Gromov-Witten invariants of $X$. For projective spaces, these invariants can be inductively computed from the Kontsevich-Manin reconstruction theorem (\cite {KM}). 

We will consider the next easiest case, that of the Grassmannians of dimension $2$ subspaces. It is well known that the genus $0$ invariants determine all higher genera (descendant) invariants. Moreover, the intersection numbers on the moduli spaces of stable bundles over curves of arbitrary genus, or certain intersection numbers on the Quot scheme can be expressed in terms of the genus $0$ Gromov-Witten invariants of Grassmannians (\cite {Ma}).

Unfortunately, Kontsevich-Manin reconstruction does not apply in this case. However, as a corollary of our study of the tautological rings, we will obtain a reconstruction result which we now explain. We begin by finding a relation between codimension $2$ evaluation classes, which is identical to the one exhibited in \cite {LP} in codimension $1$.

\begin {proposition} \label {evala} For all codimension $2$ classes $\alpha$ on a general Grassmannian $\bg$, the following codimension $2$ class on $\overline M_{0,n}(\bg, d)$: $$ev_i^{\star}\alpha-ev_j^{\star} \alpha - \psi_j \kappa(\alpha)$$ can be written explicitly as a sum of classes supported on the boundary. 
\end {proposition}

The proof of this proposition relies on a dimension computation of the complex codimension $2$ cohomology, which combines localization and the Deligne spectral sequence. As a second step, we enumerate all codimension $2$ classes on the moduli space  of maps using the main theorem in \cite {O1}. We compare the number of generators with the actual dimension of their span to explain the existence of one relation between various evaluation classes. To identify the relation explicitly we will investigate how it restricts to the space of maps to projective spaces. This will require a better understanding of the codimension $2$ classes on the space of maps to $\pr$. We will prove the following:

\begin {proposition}\label {codimensiona}
One can write down an explicit basis for the complex codimension $2$ cohomology of $\overline M_{0,0}(\p^r, d)$. 
Similar statements can be made for $1$ or $2$ marked points and for general Grassmannians. These bases are described in propositions 3.A-D. 
\end {proposition}

Once proposition $\ref {evala}$ is established, we can derive a reconstruction theorem for the genus $0$ primary Gromov-Witten invariants of the Grassmannian $\bg$ of $2$ dimensional subspaces. 

\begin {theorem}\label{reconstructiona}
(i) All genus $0$ Gromov Witten invariants of the Grassmannian of $2$ planes can be reconstructed from the invariants $\langle \alpha_1, \alpha_2, c_2, \ldots, c_2\rangle$ where $c_2$ is the second Chern class of the tautological quotient bundle and $\alpha_1, \alpha_2$ are arbitrary cohomology classes. 

(ii) For any number of marked points $n$ and any degree $d$, there is a constant $c(n,d)$ such that if $\dim V > c$ then all genus $0$, degree $d$, $n$ point Gromov Witten invariants of $\bg$ can be explicitly computed. 
\end {theorem}

We can rephrase the above result in terms of the Gromov-Witten potential as follows. Fix $\Delta_0, \Delta_1, \ldots, \Delta_l$ a basis for the cohomology of $\bg$ which is coupled to the coordinates $x_0, \ldots, x_l$ of the small phase space. We make the convention that $\Delta_0=c_2$ and we write ${\bf X}$ for the coordinates $(x_1, \ldots, x_l)$. The Gromov-Witten potential $\Phi$ is defined as the power series: $$\Phi(x)=\sum_{n}\frac{1}{n!} \langle \Delta_{i_1}, \ldots, \Delta_{i_n}\rangle x_{i_1}\cdot \ldots \cdot x_{i_n}.$$ 
\begin {corollary} The genus $0$ Gromov-Witten potential on the small phase space of $\bg$ is determined by the initial conditions: $$\Phi|_{\bf X=0},\; \partial_i\Phi|_{\bf X=0},\;\text { and }\partial_{ij} \Phi|_{\bf X=0}\; \text {for all } 1\leq i,j\leq l.$$
\end {corollary}

Similar statements can be made for a larger class of smooth projective varieties (see remark $3$). 

Of course, the genus $0$ Gromov-Witten potential satisfies complicated differential equations such as the WDVV equations. These differential equations can essentially be derived from the topological recursion relations (TRR). It would be interesting to understand how the above corollary fits in with the constraints on the potential enumerated above. We believe that corollary $1$ and proposition $\ref {evala}$ are again consequences of the TRR. Even more generally, one could perhaps hope that all relations in the tautological rings are obtained as consequences of the TRR. Theorem 1 above provides supporting evidence.

Our paper is organized as follows. In the next section, we will prove proposition $\ref {gmodpa}$ and check our result explicitly in the case of Grassmannians. In the second section we discuss the case of $SL$ flag varieties. There, we prove theorem $1$ stated above. Finally, the last part of the paper is devoted to the proof of theorem $\ref {reconstructiona}$ which will follow quite formally after establishing proposition $\ref {evala}$.  The proofs of propositions $\ref {evala}$ and $\ref {codimensiona}$ will be achieved in the third section.

{\bf Acknowledgments.} We would like to thank Alina Marian for helpful conversations. This work would have not been possible without Alina's constant encouragement. We acknowledge professor Johan de Jong who explained to us an essential aspect of the Deligne spectral sequence, which was instrumental in our computations. We gratefully thank professor Rahul Pandharipande for the interest shown in this work. 

\section {The generators for the codimension one Chow group.}

\subsection {Preliminaries.} In this subsection, we review the relevant facts about homogeneous spaces and their Schubert stratification which we will use in this paper. 

To set the stage, let $X$ be the algebraic homogeneous space $G/P$ where $G$ is a semisimple group and $P$ is a parabolic subgroup. We pick $T$ a maximal torus, $B$ a Borel subgroup such that $T\subset B \subset P \subset G$. Our convention is that the Lie algebra $\mathfrak b$ contains all the negative roots with respect to some choice of a Weyl chamber.  We let $U^{+}$ denote the unipotent subgroup of $G$ whose Lie algebra is the sum of all positive root spaces. We let $W$ be the Weyl group and $W^{\mathfrak p}$ be the Hesse diagram of $\mathfrak p$. 

We consider the decomposition of $X$ coming from the maximal torus action on $X$. This action has isolated fixed points indexed by the elements of the Hasse diagram $W^{\mathfrak p}$. The corresponding {\it plus} Bialynicki-Birula stratification coincides with the more familiar Schubert decomposition: $$G/P=\cup_{w\in W^{\mathfrak p}} U^+ \cdot (wP).$$ We let $X_w$ be the orbit $U^{+}\cdot (wP)$ and we let $Y_w$ be its closure. $Y_w$ is a subvariety of $X$ whose codimension equals the length $l(w)$ of $w$. $Y_w$ can be written as union of lower dimensional strata (where $\geq$ refers to the Bruhat ordering): $$Y_w=\cup_{w'\in W^{\mathfrak p}, w'\geq w} X_{w'}.$$ 

The codimension $1$ cells $Y_w$ are important to us. They are in one to one correspondence with the {\it simple} roots $\alpha$ of $\mathfrak g$ not contained in $\mathfrak p$, the corresponding $w$ being the reflection $s_{\alpha}$ across the wall $\alpha$. The cycle $X_{\alpha}$ corresponds to the points of $X$ which flow to $q_{\alpha}=s_{\alpha} P$ as $t\to 0$. The $Y_{\alpha}$'s can also be described by the zeros of holomorphic sections of some very ample line bundles $L_{\alpha}$. We can write any class $\beta \in H^{2}(X, \mathbb Z)$ in the form $\beta =\sum_{\alpha} d_{\alpha} \beta_{\alpha},$ where $\beta_{\alpha}$ is the codimension $1$ class $[Y_{\alpha}]=c_1(L_{\alpha})$ for each simple root $\alpha$ not in $\mathfrak p$.  

In addition, each simple root $\alpha$ determines a rational curve in $X$ joining $P$ to $q_{\alpha}=s_{\alpha} P$. The class of this rational curve is dual to the class of $Y_{\alpha}$. The rational curve can be parametrized as: $$t\to exp(tv)P, \;\text {where } v \text { is a vector in the root space of } \alpha.$$ Similarly, such $T$ invariant curves in $X$ can be generated for all positive (not necessarily simple) roots $\alpha$ which are not in $\mathfrak p$; all $T$ invariant curves passing through $P$ are obtained this way. More generally, a $T$-invariant rational curve joining two general fixed points $wP$ and $w'P$ exists provided $w'=ws_{\alpha}$ for some root $\alpha$ not in $\mathfrak p$. 

\subsection {Generators for the codimension one Chow group.}

The purpose of this subsection is to prove the following result about the generators of the rational codimension $1$ Chow group of $\fscheme$. 

\addtocounter{proposition}{-3}
\begin {proposition}\label {gmodp}

The codimension $1$ Chow group with rational coefficients of $\fscheme$ is generated by the following classes:
\begin {itemize}
\item boundary divisors of nodal maps, the degrees and marked points being distributed arbitrarily on the two components. 
\item evaluation classes $ev_{i}^{\star} c_1(L_{\alpha})$ for each $1\leq i \leq n$ and each simple root $\alpha$ which is not contained in $\mathfrak p$, 
\item kappa classes $\kappa(Y_w)$ where $l(w)=2$ (so that $Y_w$ has complex codimension $2$ in $X$). 
\end {itemize}
\end {proposition}

{\bf Proof.} As a corollary of localization and the rationality of the moduli space of stable maps we easily see that the rational Picard group, the rational codimension $1$ Chow group and the complex codimension $1$ rational cohomology all coincide, essentially because the same is true for the fixed loci. This was explained in \cite {O2} is greater detail. In this paper we will use Chow groups/cohomology interchangeably, but in the proof of this proposition it is more convenient to make use of the Chow group of $\overline M =\fscheme$. 

We let $Y$ be the complement in $X$ of all classes $Y_w$ with $l(w)\geq 2$. We also let $U$ be the open dense cell of the Schubert decomposition. It can be defined as follows: $$U=\{x\in X\; \text {such that } \lim_{t\to 0} t\cdot x \to P\}.$$ 

We consider the following subschemes of $\overline M$: 
\begin{enumerate}
\item[(a)] The codimension $1$ boundary divisors.
\item[(b)] The subscheme of maps intersecting $Y_w$ for all $w$ with $l(w)=2$. 
\item[(c)] The subscheme of maps with markings in $Y_{\alpha}$ for all simple roots $\alpha$.
\item[(d)] The subscheme of maps which cut $Y_{\alpha}$ with multiplicity higher than $1$. 
\end {enumerate}

It is clear that the complement of all these subschemes is the locus $\mathcal X$ of maps $f:(C, x_1, \ldots, x_n) \to X$ with the following properties:
\begin {itemize}
\item The domain curve is irreducible,
\item the image of the map is contained in $Y$,
\item the markings of $f$ map to $U$,
\item $f$ intersects $Y_{\alpha}$ transversally. 
\end {itemize}

We claim that the Chow group $A^{1}(\mathcal X)=0$. It follows then that the four types of classes $(a)-(d)$ span the space of divisors. To complete the proof it remains to show that the classes in item $(d)$ are indeed among the generators we enumerated in the proposition. 

The subscheme in item $(d)$ can be described as the image of the cycle  ${\overline M}_{(0,\ldots, 0, 2)}^{Y_{\alpha}}(X, \beta)$ on $\overline M_{0, n+1}(X, \beta)$ under the map $\pi$ forgetting the last marking. Here, ${\overline M}_{(0,\ldots, 0, 2)}^{Y_{\alpha}}(X, \beta)$ is the Gathmann space of stable maps with contact order at least $2$ with the very ample hypersurface $Y_{\alpha}$ (see \cite {G} for the relevant definitions). Due to the fact that $Y_{\alpha}$ is a very ample, there is an embedding $\phi: (X, Y_{\alpha})\to (\p^N, \p^{N-1})$. The following equation: $$\pi_{\star}\left[{\overline M}_{\bullet}^{Y_{\alpha}}(X, \beta)\right]=\phi^{\star}\pi_{\star} \left[{\overline M}_{\bullet}^{\p^{N-1}}(\p^{N}, \phi_{\star}\beta)\right]$$ holds in the Chow group of $\overline M_{0, n}(\p^N, \phi_{\star}\beta)$ (see \cite {G}, theorem 2.6). Here ${\overline M}_{\bullet}^{\p^{N-1}}(\p^{N}, \phi_{\star}\beta)$ denotes the corresponding Gathmann space of maps to $\p^N$ with contact order $2$ at the hyperplane $\p^{N-1}$. Therefore, it suffices to show that on $\overline M_{0, n}(\p^N, \phi_{\star}\beta)$ the pushforward class $\pi_{\star} \left[{\overline M}_{\bullet}^{\p^{N-1}}(\p^{N}, \phi_{\star}\beta)\right]$ is in the span of the corresponding classes we claimed as generators. That is, we need to show this class is a sum boundary divisors, $\kappa$ classes and evaluation classes $ev^{\star}{\mathcal O}_{\p^N}(1)$, and then observe that these classes pullback to similar classes under $\phi$. However, this statement is already proved by Pandharipande (\cite {divisors}) who in fact enumerated all divisor classes for $\overline M_{0, n}(\p^N, \phi_{\star}\beta).$

To prove the vanishing of $A^{1}(\mathcal X)$ we first consider the case when $n+\sum_{\alpha} d_{\alpha}\geq 4$. We let  $\mathfrak S=\times_{\alpha} S_{d_{\alpha}}$ and $$\mathcal F= M_{0, n+\sum_{\alpha} d_{\alpha}}/\mathfrak S, \; \overline {\mathcal F}= {\overline M}_{0, n+\sum_{\alpha} d_{\alpha}}/\mathfrak S.$$ In fact $\overline {\mathcal F}$ is the big fixed locus for the $T$ action on $\overline M$. The embedding $\overline {j}:\overline {\mathcal F} \to \overline M$ (and similarly $j:\mathcal F\to \overline M$) is obtained as follows:

\begin {itemize}
\item We consider a stable curve with $n+\sum_{\alpha} d_{\alpha}$ marked points. This will be a contracted component of the stable map whose image is the origin $P$ of $X$. We make the first $n$ marked points of the stable curve be the marked points of the stable map to $X$. 
\item At the $d_\alpha$ marked points we add $\mathbb P^1$'s of degree $1$ mapping to the rational curve joining $P$ to $q_{\alpha}$ constructed in the beginning of this section. 
\end {itemize}

We let $${\mathcal E} =\{f \text { stable map in }\overline M \text { such that } t\cdot f \to F \in  {\mathcal F} \text { as } t\to 0\}.$$ Proposition 2 of \cite {KP} shows that $\mathcal X$ is an open subvariety of $\mathcal E$. It is enough to show $A^{1}(\mathcal E)=0$. 

Let $\pi:\widehat M \to \overline M$ be a $T$-equivariant resolution of singularities for $\overline M$. The image of restricted map $j:\mathcal F \to \overline M$ lies in the smooth (automorphism-free) locus of $\overline M$. Since $\pi$ is an isomorphism over the smooth locus, we obtain an inclusion $\hat j: \mathcal F \to \widehat M$. 

Let $\widehat {\mathcal E}$ be the subset of $\widehat M$ of points flowing to $\mathcal F$. Since, $\pi$ is an isomorphism on $\mathcal F$, we have $\pi^{-1}\mathcal E = \widehat {\mathcal E}$. Now, $\pi: \widehat {\mathcal E} \to \mathcal E$ can be chosen to be a composition of blowups, so we conclude that $A^{1}(\widehat {\mathcal E})\to A^{1}(\mathcal E)$ is surjective. It is enough to show $A^{1}(\widehat {\mathcal E})=0$. This follows easily, since $\widehat M$ is smooth and for smooth varieties, it is well known that $\widehat {\mathcal E}$ is a bundle over the fixed set $\mathcal F$. Thus the claimed vanishing of $A^{1}(\widehat {\mathcal E})$ follows from the vanishing of $A^{1}(\mathcal F)=A^{1}(M_{0,n+\sum_{\alpha} d_{\alpha}})^{\mathfrak S}$. This is well known, it can be derived, for example from Keel's result that the boundary classes generate all codimension one classes on $\overline M_{0, n+\sum_{\alpha} d_{\alpha}}$.

To finish the proof we have to analyze each of the remaining cases when $n+\sum_{\alpha} d_{\alpha}\leq 3$ individually. Then $\mathcal F$ is to be interpreted as a point, and as long as this point has no automorphisms in the moduli space $\overline M$ we are done by the same arguments as before. There are four cases to consider. We will briefly show the argument for $n=0, \beta = 3\beta_{\alpha}$. The remaining cases can be obtained as in theorem 3 in \cite {KP}, by adding more marked points to place ourselves in the case we already discussed. 

We let $\mathcal Y$ be the open subscheme of $M_{0,3}(X, 3 \beta_{\alpha})$ consisting in maps with image in $Y$ such that all $3$ markings map to $Y_{\alpha}$ with multiplicity $1$. It follows that $\mathcal Y/S_3=\mathcal X$. It is therefore enough to show $A^{1}(\mathcal Y)=0$. Now, let $\widehat {\mathcal Y}$ be an equivariant resolution of singularities for the closure $\overline {\mathcal Y}$ of $\mathcal Y$ in $\overline {M}_{0,3}(X, 3 \beta_{\alpha})$.  Since $\mathcal Y$ is smooth we can view it as a subscheme of $\widehat {\mathcal Y}$. The arguments of proposition $2$ in \cite {KP} show that all $f\in \mathcal Y$ flow to a unique map $\mu \in \overline {\mathcal Y}$. This map $\mu$ has an internal component of degree $0$ mapping to $P$ to which we attach external components of degree $\beta_{\alpha}$, each of them having a marked point mapping to $q_{\alpha}$. Therefore $\mu$ sits in the smooth part of $\overline {\mathcal Y}$ and we can therefore regard it as an element in $\widehat {\mathcal Y}$. Let us look at the subvariety $\mathcal V$ of points of $\widehat {\mathcal Y}$ flowing to $\mu$ under the $\cs$ action. $\mathcal Y$ is contained in $\mathcal V$ by the above discussion. Now since $\widehat {\mathcal Y}$ is smooth, $\mathcal V$ is an affine space. Therefore $A^{1}(\mathcal Y)=0$, as desired. This completes the proof of the proposition. 

\subsection {Divisors on the space of maps to Grassmannians.} This subsection will contain another way of finding divisor classes on the space of maps to the Grassmannian $X$ parameterizing $k$ dimensional projective subspaces of $\p^r$. The method is quite ad-hoc and it involves decomposing the space of stable maps into pieces we understand better. The author could not make this procedure work for other flag manifolds. We seek to reprove the generation result of the previous subsection. 

Let $\s$ and $\q$ denote the tautological and quotient bundles on $X$. We let $W$ be a copy of $\p^{r-k-2}$ (given by the vanishing of the first $k+2$ homogeneous coordinates). Then the subvariety $$\mathcal V=\{\Lambda \in X \text { such that } \Lambda \cap W \neq \emptyset\}$$ has codimension $4$. It is easy to observe that $c_1(\q)^2$ and $\left[\mathcal V\right]$ generate the complex codimension $2$ classes on $X$. 

The complement of $\mathcal V$ parametrizes subspaces $\Lambda$ in $\p^r\setminus \p^{r-k-2}$. The key observation is that $\p^r\setminus \p^{r-k-2}$ can be understood as the total space of the bundle $$\pi: \mathcal O_{\p^{k+1}}(1)^{\oplus (r-k-1)}\to \p^{k+1}.$$ Since the fibers of $\pi$ are affine spaces, the $k$ dimensional subspace $\Lambda$ contained in $X\setminus \mathcal V$ projects to a $k$ dimensional subspace $L=\pi(\Lambda)$ in $\p^{k+1}$ i.e. $L$ gives an element of the projective space $\p^{k+1}$ of $k$ dimensional subspaces of $\p^{k+1}$. We conclude that $X\setminus \mathcal V$ can be described as a bundle $\mathcal E$ over $\p^{k+1}$ whose fiber over a subspace $L$ is $H^{0}(L, \mathcal O_{L}(1))^{\oplus (r-k-1)}$. 

An easy argument, involving the Euler sequence identifies this bundle with $T\p^{k+1}(-1)$. The long exact sequence in cohomology induced by the Euler sequence shows that this bundle is convex. That is, for any morphism $g:\p^1\to \p^{k+1}$ we have $H^{1}(\p^1, g^{\star}\mathcal E)=0$. 

We define $\mathcal X$ to be the subscheme of ${\overline M}_{0,n}(X, d)$ parameterizing stable maps $f$ whose images intersect $\mathcal V$. The open set ${\overline M}_{0,n}(X, d)\setminus \mathcal X$ consists in maps whose images are contained in the total space $X\setminus \mathcal V$ of the bundle $\mathcal E\to \p^{k+1}$. Proposition 2.1 in \cite {BH} and the convexity of the bundle $\mathcal E$ imply that ${\overline M}_{0,n}(X, d)\setminus \mathcal X$ is in fact a bundle over $\overline M_{0,n}(\p^{k+1}, d)$. Therefore, $A^{1}({\overline M}_{0,n}(X, d) \setminus \mathcal X)=A^{1}(\overline M_{0,n}(\p^{k+1}, d))$. This last group is well known, it has been computed by Pandharipande in \cite {divisors}. For example, in the case when $n\geq 3$ or $n=0$, the generators are the boundary divisors and the class $\pi_{\star} (ev_{n+1}^{\star} c_1(\mathcal O_{\p^{k+1}}(1))^2)$. This implies that the boundary divisors and the class $\pi_{\star}ev_{n+1}^{\star}c_1(Q)^2$ on $\fscheme$ restrict to generators for $A^{1}({\overline M}_{0,n}(X, d) \setminus \mathcal X)$. To get all divisors classes on $A^{1}(\fscheme)$ we need to add the class $\left[\mathcal X\right]=\pi_{\star}ev_{n+1}^{\star} \left[\mathcal V\right]$. A similar argument works for $n=1$ or $n=2$. We recover the statement of proposition $1$.

\section {The classes on the moduli spaces of maps to SL flags.}

In this section we will focus explicitly on the case of divisor classes on the space of maps to $SL$ flag varieties.  We will start by restating the results of the previous section for $SL$ flags, then describe relations between the generators we found and finally prove their independence by a dimension computation. 

\subsection {Divisors on the moduli spaces of maps to $SL_n$ flags.} To set the stage, we let $X$ be the flag variety parameterizing quotients of a vector space $V$ of fixed dimensions $n_1, \ldots, n_l$, or equivalently of subspaces of $V$ of dimensions $m_1, \ldots, m_l$: $$0\to S_1 \to \ldots \to S_l \to V \to Q_1 \to \ldots \to Q_l \to 0.$$ There is a tautological sequence on $X$ given by \begin {equation}\label {tautflag}0\to \mathcal S_1 \to \ldots \to \mathcal S_l \to V \otimes \mathcal O_{X} \to \mathcal Q_1 \to \ldots \to \mathcal Q_l \to 0.\end {equation} We let $\mathcal K_j$ denote the kernel of the map $\q_j\to \q_{j+1}$ for all $0\leq j\leq l$, where by convention $\q_0=V\otimes \mathcal O_X$ and $\q_{l+1}=0$. 

It is well known that the Chern classes $c_1(\q_j)$ form a basis for $H^{2}(X, \mathbb Z)$ for $1\leq j\leq l$. We let $h^{2}(X)$ denote the dimension of this vector space. Each stable map to $X$ will have a multi-degree $(d_1, \ldots, d_l)$ determined by the above generators of $H^{2}(X)$. Similarly,  $H^{4}(X, \mathbb Z)$ is generated by the classes $c_1(\q_i)c_1(\q_j)$ together with the {\it nonzero} Chern classes $c_2(\mathcal K_j)$. There is only one relation between these generators: \begin {equation}\label {flagcoh}\sum_{i} c_2(\mathcal K_i) - \sum_{i} c_1(\mathcal Q_i)^2 + \sum_{i} c_1(\mathcal Q_i)c_1(\mathcal Q_{i+1})=0.\end {equation} 

We enumerate the generators we obtained in proposition $\ref {gmodp}$:

\begin {enumerate}
\item [(a)] boundary classes. We have $[2^{n-1}(d_1+1)\ldots (d_l+1)]^{+}-1 - n$ such boundaries. 
\item [(b)] $\kappa$ classes $\kappa(c_1(\q_i)\cdot c_1(\q_j))$ and all classes $\kappa(c_2(\mathcal K_i))$ when $\mathcal K_i$ has rank at least $2$. 
\item [(c)] evaluation classes $ev_{i}^{\star} c_1(\q_j)$, for each $1\leq i\leq n$ and $1\leq j\leq l$.
\end {enumerate}

Now, these generators turn out not to be independent. We exhibit relations between them. We start with the boundary classes. The obvious way of getting relations is to pull back relations from $\monbar$ under the forgetful map: $$\fscheme \to \monbar.$$ There are $2^{n-1}-1-n$ boundary classes on $\monbar$ but there are $\frac{n(n-3)}{2}$ independent relations between them. This cuts down the number of independent classes in $(a)$ to at most $[2^{n-1}(d_1+1)\ldots (d_l+1)]^{+}-1-\binom{n}{2}$, with equality when all relations come from $\monbar$.

Next, we consider the classes of type $(c)$. When $n\geq 1$, we pick $H$ an ample generator and pick $m$ large enough such that the bundles ${\bf Q_j}=\det \q_j (mH)$ are all very ample. We will replace the bundles $\q_j$ in $(c)$ by their very ample counterparts ${\bf Q_j}$. Note that the span of the classes in $(b)$ and $(c)$ will be not be affected by this change. 

When $n\geq 3$ and when $X=\pr$ all divisor classes, including the corresponding evaluation classes of type $(c)$, are spanned by boundaries and the $\kappa$ class $\kappa(c_1(\mathcal O_{\pr}(1))^2)$. This is the contents of lemma 1.1.1 in \cite {divisors}. Using the linear system $|{\bf Q_j}|$ we get an {\it embedding} of $X$ into a projective space. Pulling back under this map, we can therefore conclude that the class $ev_{i}^{\star} c_1({\bf Q_j})$ is in the span of boundaries and of $\kappa$ classes, which will necessarily be on the list $(b)$. When $n\geq 3$ we will henceforth dispense with the classes $(c)$. 

When $n=1$ a different discussion is needed. We use Lemma 2.2.2 in \cite {divisors} quoted as equation $\eqref{strange}$ below. Pull back the relation provided by the lemma under the embedding given by the linear system $|{\bf Q_j}|$. Then, modulo boundary classes the following relation holds for some constants $D_j$: 
$$
\psi_1 = \frac{1}{D_j^2} \kappa(c_1({\bf Q_j})^2) - \frac{2}{D_j} ev_1^{\star} c_1({\bf Q_j}) 
$$
It follows that the span of the evaluation classes $ev^{\star}c_1({\bf Q_j})$ is exactly $1$ dimensional modulo boundaries and the $\kappa$ classes in $(b)$ (for example $\psi_1$ is a generator of the one dimensional span. Intersecting with suitable curves or restricting to a copy of $\p^1\hookrightarrow X$ and invoking the results of \cite {divisors}, it can be shown this class is independent from the boundaries). 

For $n=2$, a similar discussion as above shows that the span of the evaluation classes in $(c)$ is at most $2$ dimensional modulo boundaries and the kappa classes in $(b)$. The classes $\psi_1$ and $\psi_2$ can be chosen as generators for the two dimensional span (alternatively we can pick pairs of evaluation classes). However, corollary 1 in \cite {LP} rewritten as equation $\eqref {psisum}$ below shows that the sum of these two classes is also in the span of the boundaries. It follows that the span of the evaluation classes $(c)$ is exactly $1$ dimensional modulo boundaries and $\kappa$ classes. We immediately conclude that for all values of $n$, the combined contribution of the classes in $(a)$ and $(c)$ is at most $[2^{n-1}(d_1+1)\ldots (d_l+1)]^{+}-1-\binom{n}{2}$. 

Finally, the classes of type $(b)$ are also connected by relations. Lemma $\ref {relation}$ of the next subsection shows that $\kappa(c_1(\q_i) c_1(\q_j))$ can be expressed in terms of $\kappa(c_1(\q_i)^2)$ and $\kappa(c_1(\q_j)^2)$ modulo boundaries. In fact, applying this lemma to each pair $(\mathcal Q_{i}, \mathcal Q_{i+1})$ we derive that the following equation is true modulo boundaries: $$\kappa(c_1(\mathcal Q_i) c_1(\mathcal Q_{i+1}))=\frac{d_{i+1}}{2d_i} \kappa(c_1(\mathcal Q_i)^2)+\frac{d_i}{2d_{i+1}} \kappa (c_1(\mathcal Q_{i+1})^2).$$ Using $\eqref {flagcoh}$, we easily arrive at the following relation: \begin {equation}\label {flageq} \sum_{i} \kappa(c_2(\mathcal K_i))+ \sum_{i} \left(\frac{d_{i-1}+d_{i+1}}{2d_i}-1\right) \kappa(c_1(\mathcal Q_i)^2) =0 \text { modulo boundaries}. \end {equation}
The coefficients of the boundary terms can be written down explicitly, but they

To summarize, we obtain the following upper bound for the dimension of the space of divisors on $\fscheme$: 
\begin {equation}\label {upper}
\left [2^{n-1}(d_1+1) \ldots (d_l+1)\right]^{+} - 1 - \binom {n}{2}+h^{4}(X) - \binom {h^2(X)}{2}.
\end {equation}

In section $2.4$, we will use localization to give a lower bound for the dimension of the space of divisors. In fact, we will prove that the bound obtained above is sharp. This will give the proof of theorem $1$ stated in the introduction.

\subsection {Relations between the $\kappa$ classes.} We will now indicate the statement and proof of the lemma invoked in the  previous subsection to find relations between the $\kappa$ classes. We hope this lemma could also be of use to understand divisors on product spaces. 

\begin {lemma} \label {relation} Let $L, M$ be two line bundles on a projective variety $X$. Then the following $\kappa$ class on ${\overline M}_{0,0}(X, \beta)$: $$\kappa\left(\left(\frac{c_1(L)}{\int_{\beta} c_1(L) } - \frac{c_1(M)} {\int_{\beta} c_1(M)} \right)^2\right)$$ is in the span of the boundary divisors. \end {lemma}

{\bf Proof.} We claim that it is enough to prove the statement for $L$ and $M$ very ample. Indeed, assuming we proved the statement in this case, we pick a very ample divisor $H$ on $X$. For $n$ large enough, $L+nH, M+nH$ will both be very ample. We obtain that $$\kappa\left(\left(\frac{c_1(L)+nH}{\int_{\beta} c_1(L) + n H\cdot\beta}-\frac{c_1(M)+nH}{\int_{\beta} c_1(M) + n H\cdot\beta}\right)^2\right)$$ is in the span of boundaries. Clearing denominators, and then looking at the term independent of $n$, we derive that the $\kappa$ class in the statement of the lemma is also in the span of boundaries. 

Assume now $L$ and $M$ are very ample. We consider the embedding $i:X \to \mathbb P^n \times \mathbb P^m$ determined by the linear systems $|L|$ and $|M|$. We let $d=\int_{\beta}c_1(L)$ and $e=\int_{\beta} c_1(M)$. We let $\mathcal H_1$ and $\mathcal H_2$ be the two hyperplane bundles on the projective spaces $\mathbb P^n$ and $\mathbb P^m$. Let $D_{i,j}$ denote the boundary divisor of maps with nodal target such that the bidgree of the map on one of the components is $(i,j)$; then the bidgree on the other component is $(d-i, e-j)$. 

The lemma will then follow pulling back under $i$ the following relation on $\mpnpmbar$: 
\begin {equation}\label {reconstruction}
\kappa \left(\left(\frac{c_1(\mathcal H_1)}{d} - \frac {c_1(\mathcal H_2)}{e}\right)^2\right)
 = \frac{1}{2} \sum_{i=1}^{d}\sum_{j=1}^{e} D_{i,j} \left(\frac{i}{d}-\frac{j}{e}\right)^2
\end {equation} 
{\bf Claim.} As a first step in establishing $\eqref {reconstruction}$, we show that the codimension $1$ classes on ${\overline M}_{0,0}(\p^n\times \p^m, (d,e))$ are in the span of the boundaries and of the two kappa classes $\kappa (c_1(\mathcal H_1)^2)$ and $\kappa (c_2(\mathcal H_2)^2)$.

The boundary in $\mpnpmbar$ is a divisor with normal crossings. It follows from the Deligne spectral sequence that the cokernel of the Gysin map $$\oplus H^{0}(\text {boundaries})\to H^{2}(\mpnpmbar)$$ can be identified with the weight $2$ piece of the Hodge structure on the cohomology of the open stratum $W^{2}H^{2}(\mpnpm)$. We will show this is at most two dimensional. 

We will write the open stratum $M=\mpnpm$ as a global quotient and make the computation in equivariant cohomology. 
Let $V$ be a two dimensional space so that $\p^1=\p(V)$ comes with the obvious $PGL(V)$ action. The space of maps $\text {Map}=\text {Map}_{(d,e)}(\p^1, \p^n\times \p^m)$ of bidgree $(d,e)$ is an open set in the product of two projective spaces $$\p\left(\bigoplus_{i=0}^{n} {Sym}^{d} V^{\star} \right)\times \p \left(\bigoplus_{i=0}^{m} {Sym}^{e} V^{\star} \right).$$ We need to factor out the action of $PGL(V)$ on the two factors to obtain $M_{0,0}(\p^n\times \p^m, (d,e))$. Equivalently, we can think of $\mpnpm$ as sitting in a quotient of the affine space $$\mathbb A = \bigoplus_{i=0}^{n} {Sym}^{d} V^{\star} \oplus \bigoplus_{i=0}^{m} {Sym}^{e} V^{\star}$$ by the action of the group $GL (V) \times \cs$. The action of $GL (V)$ is the usual one on the two factors, while $\cs$ acts in the usual way only on the second factor, and trivially on the first. This action is easily seen to have finite stabilizers. It is well known that in such cases we have an isomorphism between the cohomology of the orbit space and equivariant cohomology: $$H^{\star}(\mpnpm)=H^{\star}(\text {Map}\times_{GL_2 \times \cs} (EGL_2 \times E\cs)).$$ Both sides have Hodge structures (for equivariant cohomology, we need to use finite dimensional approximations of the equivariant models) compatible with the above isomorphism. Moreover, $\text {Map}\times_{GL_2 \times \cs} (EGL_2 \times E\cs)$ sits inside the space $$\left(\bigoplus_{i=0}^{n} {Sym}^{d} V^{\star} \oplus \bigoplus_{i=0}^{m} {Sym}^{e} V^{\star} \right)\times_{GL_2 \times \cs} (EGL_2 \times E\cs).$$ This space is the total space of a bundle over the product of two classifying spaces $BGL_2 \times B\cs$. The restriction map $$H^{2}(BGL_2\times B\cs)=H^{2}_{GL_2\times \cs}(\mathbb A)\to W^{2}H^{2}_{GL_2\times \cs}(\text {Map})=W^{2}H^{2}(M)$$ is surjective. Therefore $W^{2}H^{2}(M)$ is at most $2$ dimensional. The surjectivity of the map can be explained by the usual arguments in \cite {O1}, using (Grothendieck's) remark $3$ of that paper. 

Our {\bf claim} follows if we show that the two classes $\kappa (c_1(\mathcal H_1)^2)$ and $\kappa (c_2(\mathcal H_2)^2)$ are not in the linear span of the boundary divisors. This is done in \cite {divisors}, Lemma 1.2.1(i) by intersecting with curves in the moduli space in the case of $\p^r$, but the argument goes through without change for $\p^n\times \p^m$. 

The claim we just proved suffices to establish the results needed in this paper. However for the sake of completeness, we will also prove the precise relation $\eqref {reconstruction}$. It follows from what we proved above that a linear combination of the three classes: \begin {equation} \label {rec} \kappa (c_1(\mathcal H_1) c_1(\mathcal H_2)) + A \cdot \kappa (c_1(\mathcal H_1)^2) +B \cdot \kappa (c_1(\mathcal H_2)^2)=\text {sum of boundary classes}\end {equation} To identify the coefficients of this relation it is enough to intersect $\eqref {rec}$ with curves in ${\overline M}_{0,0}(\mathbb P^n \times \mathbb P^m, (d,e))$. A moment's thought shows that it is enough to check $\eqref {reconstruction}$ for all curves of $\overline M$ transversal to the boundary. Indeed, assuming this is the case, we show that  $\eqref {rec}$ and an appropriately scaled version of $\eqref {reconstruction}$ coincide. Subtracting the two equations, we get an expression involving only $\kappa (c_1(\mathcal H_1)^2)$, $\kappa (c_1(\mathcal H_2)^2)$ and boundary classes. This expression vanishes on each curve in $\overline M$ transversal to the boundary. We have seen already in the proof of the claim that this implies that the coefficients of the $\kappa$'s must vanish. It is not any harder to conclude the same about the coefficients of the boundary classes. \footnote {Recall that as a consequence of localization, numerical and rational equivalence coincide, essentially because the same is true for each of the fixed loci. One can perhaps conceive an argument which would establish $\eqref {reconstruction}$ without appealing to the {\bf claim} above, simply by intersecting with the smooth curves of $\overline M$. We actually do this in the proof below for curves transversal to the boundary. The general argument should not be more complicated.} 

It remains to show $\eqref {reconstruction}$ holds after intersecting with the smooth curves intersecting the boundary divisors transversally. Let us now consider such a curve. This is the same as a family of stable maps to $\mathbb P^n \times \mathbb P^m$ parametrized by a one dimensional base $B$: 
\begin {center} 
$\begin {CD}
S @>{F=(f,g)}>>\mathbb P^n \times \mathbb P^m\\
@V{\pi}VV \\
B
\end {CD}$
\end {center} 
It can be proved that $S$ is the blow up of a projective bundle $P=\p(V)$ at the points $x_1, \ldots, x_s$ where $B$ meets the boundary divisors. 

We let $E_i$ be the exceptional divisors of the blowups and we let $h=c_1(\mathcal O_{\p(V)}(1))$. We assume that the map $F$ has bidegree $(d_i, e_i)$ on each exceptional divisor $E_i$. It is then clear that for some line bundles $\mathcal J_1$ and $\mathcal J_2$ on $B$ we have: \begin {equation}\label {pull}F^{\star} \mathcal H_1 = \pi^{\star} \mathcal J_1 \otimes \mathcal O_{\p(V)}(d) \otimes \mathcal O(-\sum_{i} d_i E_i)\end {equation} \begin{equation}\label {pull1}F^{\star} \mathcal H_2 = \pi^{\star} \mathcal J_2 \otimes \mathcal O_{\p(V)}(e) \otimes \mathcal O(-\sum_{i} e_i E_i).\end {equation} 

It is also obvious that $B\cdot D_{i,j}=n(i,j)+n(d-i, e-j)$ where $n(u,v)$ is the number of points among $x_1, \ldots, x_s$ such that the bidegree of the map $F$ on the corresponding exceptional divisor is $(u,v)$. 

We now intersect both sides of the equation ($\ref {reconstruction}$) with the curve $B$. We will need to show that \begin {eqnarray*}
\pi_{\star}\left( \left(\frac{1}{d}c_1(F^{\star}\mathcal H_1)-\frac{1}{e}c_1(F^{\star}\mathcal H_2)\right)^2\right) &=&\frac{1}{2} \sum_{i,j} (n(i,j)+n(d-i, e-j))\left(\frac{i}{d}-\frac{j}{e}\right)^2\\&=&\sum_{i,j} n(i,j)\left(\frac{i}{d}-\frac{j}{e}\right)^2=\sum_{i}\left(\frac{d_i}{d}-\frac{e_i}{e}\right)^2\end {eqnarray*} 

The proof of this equality involves equations ($\ref {pull}$) and ($\ref {pull1}$). Indeed, the right hand side of the expression above equals $$\pi_{\star} \left(\left(\frac{1}{d} (\pi^{\star} c_1(\mathcal J_1)+d h - \sum_{i} d_i E_i) - \frac{1}{e} (\pi^{\star} c_1(\mathcal J_2)+e h - \sum_{i} e_i E_i)\right)^2\right)$$ After a few cancellations, we finally arrive at the desired result. The proof of the lemma is complete. 

\subsection {The computation of the symmetric group invariants.} Our next digression will be useful in the dimension computation needed to finish the proof of theorem $\ref {diva}$. We will prove a preliminary result about the $S_n$ action on the cohomology of the moduli space of rational pointed curves $\monbar$.

To fix notation, for each permutation $\sigma \in S_n$ we write $n_j(\sigma)$ for the number of cycles of length $j$. We denote by $c(\sigma)$ the total number of cycles of $\sigma$. 

\begin {lemma}\label {h2mon}
For each $\sigma \in S_n$, the trace of $\sigma$ on $H^{2}(\monbar)$ is given by $$2^{c(\sigma)-1}-1-n_2(\sigma)-\binom{n_1(\sigma)}{2}+\delta(\sigma)$$ where $$\delta(\sigma)=\begin {cases} 2^{c(\sigma)-1} & \text {if $\sigma$ has only cycles of even length} \\ 0 & \text {otherwise} \end {cases}$$
\end {lemma}

{\bf Proof.} The proof of this lemma makes use of the ideas of Getzler's paper (\cite {zero}). Getzler works out the Deligne spectral sequence of the mixed Hodge structure on the open manifold $M_{0,n}$. He shows that there is an exact sequence:
\begin {equation}\label {getz} 0\to H^1(\mon)\to \oplus_{\Gamma} H^{0} (D_\Gamma)\to H^2(\monbar)\to 0.\end {equation}
Here $D_{\Gamma}$ are the boundary divisors of $\monbar$. They correspond to unordered partitions of the $n$ marked points into two subsets $A$ and $B$ such that $|A|, |B|\geq 2$. 

It is clear that the the middle term of the exact sequence $\eqref {getz}$ is a sum of one dimensional spaces, one for each unordered partition of $\{1, \ldots, n\}$ into $2$ subsets as above. If $n$ is odd, the trace of $\sigma$ on the middle term of the exact sequence $\eqref {getz}$ equals the number of partitions $\{A,B\}$ such that $$\sigma(A)=A,\;\sigma(B)=B,\; |A|, |B|\geq 2$$ This number is easily seen to be $2^{c(\sigma)-1}-1-n_1(\sigma)$. Indeed, both $A$ and $B$ have to be unions of full cycles of $\sigma$. The last two terms are the corrections corresponding to the non-stable cases when $A$ or $B$ have $0$ or $1$ elements. In the case when $n$ is even, we also need to consider the partitions $\{A, B\}$ such that $$\sigma(A)=B, \sigma(B)=A, n\geq 4$$ Such partitions exist only if $\sigma$ has all cycles of even length and their number equals $2^{c(\sigma)-1}$. 

The proof will be complete using the exact sequence $\eqref {getz}$, the remarks above and if we also show that:
\begin {equation}\label {trace}
\text {Tr}_{\sigma } H^1(\mon)=n_2(\sigma)+\frac{n_1(\sigma)^2-3n_1(\sigma)}{2}
\end {equation}

To establish $(\ref {trace})$, we will use the following facts collected from \cite {zero}. \\
\begin {enumerate}\item[(a)] First, Getzler shows that as a consequence of the Serre spectral sequence, we have an isomorphism $H^{\bullet}(M_{0,n+1})=H^{\bullet}(M_{0,n}\times \mathbb C\setminus\{1,2\ldots, n\})$. It is then clear that $$\text {Tr}_{\sigma }H^1(M_{0,n+1})=\text {Tr}_{\sigma }H^{1}(\mon)+\text {Tr}_{\sigma}H^{1}(\mathbb C\setminus\{1,2\ldots, n\})=\text {Tr}_{\sigma }H^{1}(\mon) + n_1(\sigma) - 1.$$
\item[(b)] Getzler also shows that if $F_n$ denotes the configuration space of $n$ pairwise distinct points in $\mathbb C$, then $H^{\bullet}(F_n)=H^{\bullet}(M_{0,n+1}\times S^1)$. Therefore, $$\text {Tr}_{\sigma} H^{1}(M_{0,n+1})=\text {Tr}_{\sigma} H^{1}(F_n)-1.$$
\item[(c)] Finally, it is a consequence of a formula of Lehrer and Solomon, also discussed in \cite {Ge2} that $\text {Tr}_{\sigma }H^{1}(F_n)=n_2(\sigma)+\binom{n_1(\sigma)}{2}$. These three items together prove equation $(\ref {trace})$, thus completing the proof of the lemma.
\end {enumerate}

\begin {lemma}\label {dim}
Let $n, a_1, \ldots a_l$ be positive integers. Consider the obvious action of $S_{a_1}\times \ldots \times S_{a_{l}}$ on $H^{2}({\overline M}_{n+a_1+\ldots+a_l})$. The dimension of the invariant subspace is computed by the formula: $$\text {dim } H^{2}({\overline M}_{n+a_1+\ldots+a_l})^{S_{a_1}\times \ldots \times S_{a_{l}}}=\left[2^{n-1}(a_1+1)\ldots (a_l+1)\right]^{+}-1-\binom {n}{2}-ln-\binom {l+1}{2} + \mathfrak a.$$ Here $\mathfrak a$ denotes the number of indices $i$ such that $a_i=1$. We also write $[x]^{+}=x$ if $x$ is an integer and $[x]^{+}=x+\frac{1}{2}$ if $x$ is a half integer. 
\end {lemma}

{\bf Proof.} To prove this statement, we will average out the trace of each permutation $\sigma \in S_{a_1} \times \ldots \times S_{a_{l}}$ on $H^{2}({\overline M}_{n+a_1+\ldots+a_l})$. For this combinatorial computation we will need the following identities which can be proved by induction on $k$.  

\begin {equation}\label {cy}
\sum_{\sigma\in S_k} 2^{c(\sigma)}=(k+1)!
\end {equation}

\begin {equation}\label {even}
\sum_{\sigma\in S_k} \delta(\sigma)=\begin {cases} \frac{k!}{2} & \text {if $k$ is even} \\ 0 &\text {otherwise}\end {cases}
\end {equation}

\begin {equation} \label {n1}
\sum_{\sigma \in S_k} n_1(\sigma)=k!
\end {equation}

\begin {equation} \label {n2}
\sum_{\sigma \in S_k} \left(n_2(\sigma)+\binom {n_1(\sigma)}{2}\right)=k! \text { for } k\geq 2.
\end {equation}

Another induction, this time on $l$, making use of equation ($\ref {n1}$) gives the following two identities:

\begin {equation}\label {n1g}
\sum_{\sigma \in S_{a_1}\times \ldots \times S_{a_l}} n_1(\sigma)=l a_1! \ldots a_l!
\end {equation}

\begin {equation} \label {n1n1}
\sum_{1\leq i<j\leq l}\sum_{\sigma_i\in S_{a_i}}\sum_{\sigma_j \in S_{a_j}} n_1(\sigma_i) n_1(\sigma_j) = \binom {l}{2} a_1! \ldots a_l!
\end {equation}

We can now compute the dimension of the invariant subspace.  A permutation $\sigma\in S_{a_1}\times \ldots \times S_{a_l}$ is tantamount to $l$ permutations $\sigma_i\in S_{a_i}$. Lemma $\ref {h2mon}$ shows that  
\begin {eqnarray*}
\text {Tr}_{\sigma } H^2({\overline M}_{n+a_1+\ldots+ a_l})= 2^{\sum_{i}c(\sigma_i)+n-1} - 1- \sum_{i=1}^{l} n_2(\sigma_i) - \binom{n+n_1(\sigma_1) + \ldots + n_1(\sigma_l)}{2} \\
+ \begin {cases} 2^{\sum_i c(\sigma_i) - 1} &\text {if $n=0$ and the $\sigma_i$'s have only even length cycles} \\ 0 &\text {otherwise} \end {cases}.
\end {eqnarray*}

We average out these traces to compute the dimension of the invariant subspace. First, we assume $n\neq 0$.

\begin {eqnarray*}
\text {dim} &H^{2}&({\overline M}_{n+a_1+\ldots+ a_l})^{S_{a_1}\times \ldots \times S_{a_l}}=
\frac{1}{a_1!\ldots a_l!} \sum_{\sigma\in S_{a_1}\times \ldots \times S_{a_l}} \text {Tr}_{\sigma} H^{2}({\overline M}_{n+a_1+\ldots+ a_l})\\&=&\frac{1}{a_1!\ldots a_l!} \left ( \sum_{i=1}^{l}\sum_{\sigma_i\in S_{a_i}} 2^{\sum_{i}c(\sigma_i)+n-1} - 1- \sum_{i=1}^{l} n_2(\sigma_i) - \binom{n+n_1(\sigma_1)+\ldots+n_1(\sigma_l)}{2}\right)\\&=&2^{n-1}\prod_{i=1}^{l} \left(\sum_{\sigma_i \in S_{a_i}} \frac{2^{c(\sigma_i)}}{a_i!}\right) -1 - \sum_{i=1}^{l}\frac{1}{a_i!}\left(n_2(\sigma_i)+\binom{n_1(\sigma_i)}{2}\right) - \binom{n}{2} - \\
&-&\frac{n}{a_1!\ldots a_l!} \sum_{i=1}^{l} \sum_{\sigma_i\in S_{a_i}} n_1(\sigma_i) - \frac{1}{a_1!\ldots a_l!} \sum_{1\leq i<j\leq l} n_1(\sigma_i)n_1(\sigma_j) =\\&=& 2^{n-1} (a_1+1)\ldots(a_l+1) - 1 - (l -\mathfrak a) - \binom {n}{2} - nl - \binom {l}{2}.
\end {eqnarray*}
where in the last line we used equations ($\ref {cy}$), ($\ref{n2}$), ($\ref{n1g}$) and ($\ref{n1n1}$) respectively.

In the case when $n=0$, we have extra-terms corresponding to the case when all permutations $\sigma_i$ have even cycles. The contribution of these terms is  
$$\frac{1}{a_1!\ldots a_l!} \sum_{\sigma_i \in S_{a_i} \text {with only even cycles}} 2^{c(\sigma_1)+\ldots +c(\sigma_l)-1}=\frac{1}{2} \prod_i \left(\sum_{\sigma_i \in S_{a_i}}\frac{2}{a_i!}\delta(\sigma_i)\right) =\frac{1}{2},$$ by virtue of equation ($\ref {even}$). This completes the proof of the lemma. 

\subsection {The fixed loci contributions.} We will exhibit particular fixed point loci for the $\cs$ action on $\fscheme$ which comes from a torus action on $X$. We will argue that the normal bundles of the fixed loci we exhibit have at most $1$ negative weight. We will employ the "homology basis theorem" to give a lower bound for the dimension of the space of codimension $1$ classes (we refer the reader to lemma $1$ in \cite {O2} for the precise statement of the homology basis theorem in the case of finite quotient singularities). The lower bound will match the upper bound $\eqref {upper}$, completing the proof of theorem $\ref {diva}$. Therefore, the argument will also give {\bf all} fixed loci which contribute to the computation of $H^{2}(\fscheme)$. 

We consider an action of $\cs$ on $V$ with weights $\lambda_1 < \ldots < \lambda_N$ and weight vectors $e_1, \ldots, e_N$. In fact we will need assume that the weight $\lambda_i$ is much bigger than $\lambda_{i-1}$. This assumption will be needed when evaluating the negative weights for the normal bundles to the fixed loci below, especially when dealing with the vertices of valency $2$. We write our flag variety as a quotient $X=SL(V)/P$ where $P$ is the parabolic of upper triangular block matrices of size $m_1$, $m_2-m_1$, \ldots, $m_{l}-m_{l-1}$, $N-m_{l}$. We will frequently use the notation $W(i)$ for the vector subspace of $V$ spanned by $e_1, \ldots, e_i$. We enumerate the following fixed points of the $\cs$ action on $X$:\\
(a) the origin $P$. The tangent space $T_{P} X$ has no negative $\cs$ weights. Explicitly, this fixed point is represented by the flag: $$W(m_1)\subset \ldots \subset W(m_l)\subset V.$$    
(b) There are $h^{2}(X)$ fixed points $q_1, \ldots, q_l$ corresponding to the simple roots not in the parabolic subalgebra. There is only one negative weight on the tangent spaces $T_{q_i}X$. Each $q_i$ can be joined to $P$ by a rational curve $R_i$ whose Poincare dual is the generator $\beta_i$ of $H^{2}(X, \mathbb Z)$. Explicitly $q_i$ is represented by a flag which differs from the one in the previous item only at the $i$ th step: $$W(m_1) \subset \ldots \subset W(m_i-1) \oplus \text { span } \langle e_{m_i+1}\rangle \subset W(m_{i+1})\subset \ldots \subset W(m_l) \subset V.$$ Additionally, the curve $R_i$ can be parametrized as: \begin {equation}\label {line}[t:s] \mapsto \text { the flag } W(m_1) \subset \ldots \subset W(m_i-1) \oplus \langle t e_{m_i} + s e_{m_i+1}\rangle \subset W(m_{i+1})\subset \ldots \subset W(m_l)\end {equation}
(c) There are $h^{4}(X)$ fixed points with $2$ negative roots on their tangent spaces. There are three types of such points which we now describe. Let $m_0=0$ and $m_{l+1}=n$. 
\begin {itemize} \item for each pair $1\leq i< j\leq l$ such that $m_j-m_i\geq 2$ we have a fixed point denoted  $q_{ij}$. It is obtained from the reference flag in (a) by modifying its $i$th and $j$th steps: $$W(m_1)\subset \ldots \subset W(m_i-1)\oplus\langle e_{m_i+1}\rangle \subset \ldots \subset W(m_j-1)\oplus\langle e_{m_j+1}\rangle \subset \ldots \subset W(m_l) \subset V$$ It is clear that $q_{ij}$ can be joined to both $q_i$ and $q_j$ by rational curves in the cohomology classes dual to $\beta_j$ and $\beta_i$. 
\item For each $1\leq i \leq l$ such that $m_{i+1}-m_{i}\geq 2$ we obtain a fixed point $r_i$ which has the property that it can be joined to $P$ by a rational curve in the cohomology class dual to $\beta_i$. The $r_i$'s can be obtained from  the reference flag in (a) by modifying its $i$th step : $$W(m_1)\subset \ldots \subset W(m_i-1)\oplus\langle e_{m_i+2}\rangle \subset W(m_{i+1})\subset \ldots \subset W(m_l)\subset V.$$ In addition, for all $1\leq i\leq l$ such that $m_i-m_{i-1}\geq 2$, we get the fixed points $r'_i$ which can be joined to the origin by a curve in the cohomology class dual to $\beta_i$. They can be defined modifying the $i$th step of the reference flag in (a):
$$W(m_1)\subset \ldots \subset W(m_i-2)\oplus \langle e_{m_i}, e_{m_i+1} \rangle \subset W(m_{i+1})\subset \ldots W(m_l)\subset V.$$ 
\item For each $m_{i+1}-m_{i}=1$, $1\leq i < l$, we have a fixed point $s_i$ which can be joined to $q_{i+1}$ by a rational curve in the cohomology class dual to $\beta_{i}$. Explicitly it is given by changing the $i$th and $(i+1)$st steps of the standard flag: $$W(m_1)\subset \ldots \subset W(m_i-1)\oplus \langle e_{m_i+2} \rangle \subset W(m_i)\oplus \langle e_{m_i+2}\rangle \subset\ldots W(m_l)\subset V.$$ Similarly, for all $1<l\leq l$ such that $m_i-m_{i-1}=1$, we obtain the fixed point $s_i'$ which can be joined to $q_{i-1}$ by a rational curve in the cohomology class dual to $\beta_i$. Explicitly, this point is obtained by modifying the $(i-1)$st and $i$th steps of the reference flag $$W(m_1)\subset \ldots \subset W(m_{i}-2)\oplus \langle e_{m_{i}}\rangle\subset W(e_{m_{i}-2}) \oplus \langle e_{m_{i}}, e_{m_{i}+1} \rangle \subset \ldots \subset W(m_l)\subset V.$$
\end {itemize}

Turning to the fixed loci on $\fscheme$, we will employ the usual method of bookkeeping the fixed loci by means of 
decorated graphs $\Gamma$. The vertices of $\Gamma$ are in one to one correspondence with the components of $f^{-1}(\cs \text { fixed points on $X$})$. The vertices come with labels corresponding to the $\cs$ fixed points on $X$. The edges of the graph correspond to non-contracted components of the map $f$ and are decorated with the degree on that component. The graph $\Gamma$ has legs attached to its vertices. A flag $\mathfrak f$ determines an edge $e(\mathfrak f)$ and thus a non contracted component $C_e$ of the stable map. We let $R_e$ be the image of this component. The flag $\mathfrak f$ also determines a unique vertex $v(\mathfrak f)$ which gives a point on the curve $C_e$ mapping to a fixed point of the $\cs$ action on $X$. We let $\omega_{\mathfrak f}$ denote the $\cs$ weight on the fiber of the bundle $f^{\star} T R_e$ at the point $v(\mathfrak f)$.

The weights on the normal bundles of each fixed locus were computed by Kontsevich-Graber-Pandharipande (see \cite {GP}). The list of weights on the normal bundle of the fixed locus indexed by $\Gamma$ is generated by the following algorithm:
\begin {itemize}
\item flag contributions: for each flag $\mathfrak f$ whose vertex $v(\mathfrak f)$ has total valency $\geq 3$ we include the weight $\omega_{\mathfrak f}$ on our list of weights.
\item vertex contributions: for each vertex $v$ corresponding to a fixed point $q$ of the $\cs$ action on $X$, we include the $\cs$ weights on $T_{q}X$. 
\item vertex contributions: the vertices $v$ with valency $2$ and no legs have two incident flags $\mathfrak f_1$ and $\mathfrak f_2$. We include the weight $\omega_{\mathfrak f_1}+\omega_{\mathfrak f_2}$.
\item edge contributions: for each edge $e$, we include the weights of the $\cs$ action on $H^{0}(C_e, f^{\star} TX)$. 
\item flag contributions: for each flag $\mathfrak f$, the vertex $v(\mathfrak f)$ maps to a fixed point $q$ of the $\cs$ action on $X$. We remove the weights on $T_{q}X$ from the list.
\item vertex contributions: for each vertex $v$ with valency $1$ and no legs, i.e. those vertices contained in only one flag $\mathfrak f$, we remove the weight $\omega_{\mathfrak f}$ from the list.
\end {itemize}

The graph below corresponds to a fixed locus with no negative weights on its normal bundle.
\begin {figure}[htb]
\scalebox{.8}{\psfig {figure=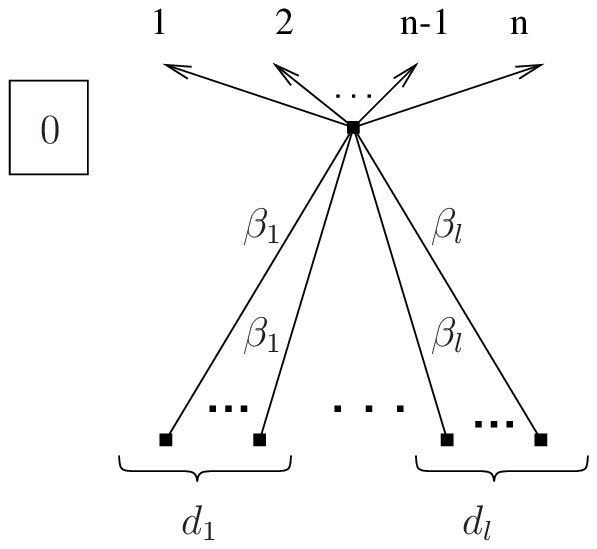}}
\caption {The big locus of the $\cs$ flow on $\fscheme$}
\end {figure}
Its {\it plus} cell is the big locus of the $\cs$ flow on $\fscheme$. The corresponding fixed locus is isomorphic to the quotient ${\overline M}_{0, n + \sum_{i} d_i}/{S_{d_1}\times \ldots \times S_{d_l}}$. Its contribution to $H^{2}(\fscheme)$, as determined by the "homology basis theorem" in \cite {O2}, equals the second Betti number. By lemma $\ref {dim}$ this contribution is: $$[2^{n-1}(d_1+1)\ldots (d_l+1)]^{+}-1-ln - \binom{l+1}{2}-\binom{n}{2}+\mathfrak a,$$ where $\mathfrak a$ is the number of indices such that $d_j=1$. 

There are (at least) five types of graphs which correspond to fixed loci with only one negative weight on their tangent bundle.  To compute their contribution to $H^{2}(\fscheme)$, we only have to count all such graphs. 
\begin {figure}[htb]
\scalebox {.8}{\psfig {figure=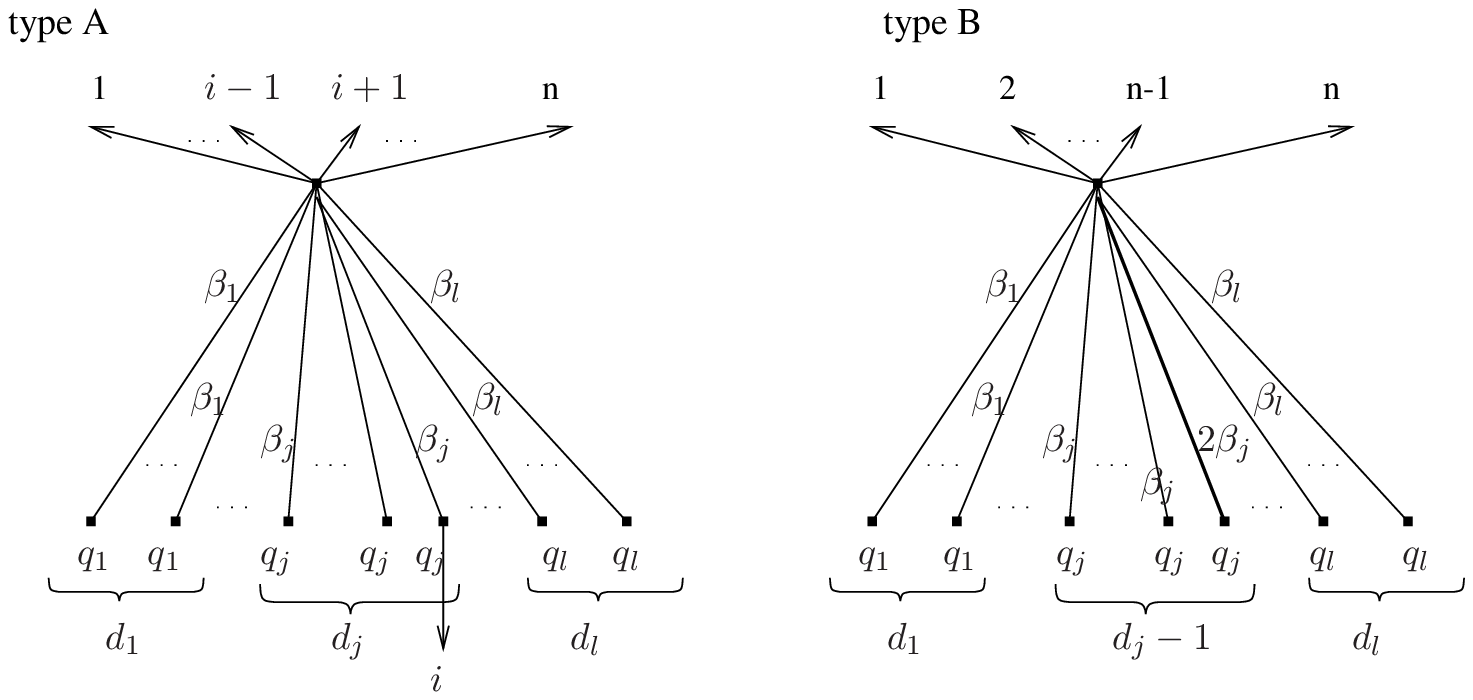}}
\scalebox {.8}{\psfig {figure=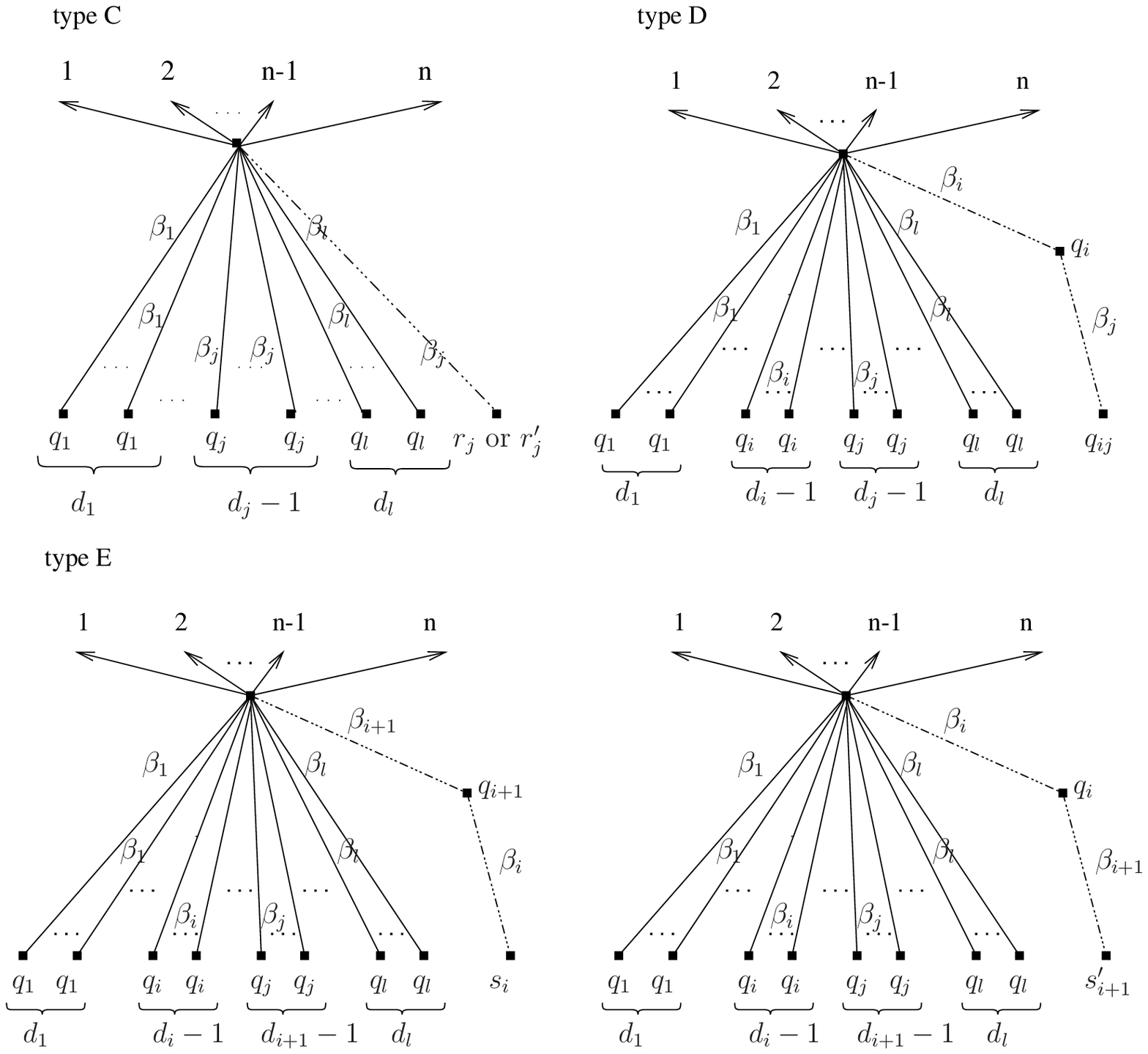}}
\caption {Fixed loci with one negative weight on the normal bundle.}
\end {figure}

We will analyze each of the five types one by one. In drawing these graphs, we used continuous lines for the edges representing the rational curves $R_i$ in the cohomology class dual to $\beta_i$, which join the origin $P$ to the fixed points $q_i$. We also indicated just below the graph the number of such curves that we use. 

The graphs of type $A$ have one leg labeled $i$ attached to a vertex labeled by $q_j$. There are $n \cdot l$ such graphs. The graphs of type $B$ have one thicker edge labeled $2\beta_j$. This edge corresponds to components mapping to $R_j$ with degree $2$ with ramification only over $P$ and $q_j$. These graphs only exist if $d_j\geq 2$, and their number equals $l-\mathfrak a$. For the graphs of type $C$, one of the edges corresponding to the rational curve $R_j$ is replaced by a rational curve, still in the homology class $\beta_j$ which joins the origin to $r_j$ or $r'_j$. This new edge is represented by a dotted line. For each $1\leq i < j \leq l$ we obtain a graph of type $D$. An edge representing the curve $R_j$ has been replaced by a rational curves with cohomology class $\beta_j$ joining $q_i$ to $q_{ij}$. The new edge is attached to one of the edges representing the curve $R_i$. We used dotted lines for the two rational curves in question. Notice the apparent asymmetry between $i,j$. Indeed, the corresponding graph obtained by switching $i$ and $j$ has one more negative weight. This comes from the contribution of the vertex of valency $2$. In our case, that vertex contributes with positive weight as one immediately checks remembering that the weight $\lambda_{j+1}$ is the dominant one among the weights which appear on the rational curves coming into the valency $2$ vertex in question. Finally, the graphs of type $E$ are described in the same manner. We have replaced two rational curves in the cohomology classes $\beta_{i}$ and $\beta_{i+1}$ by two rational curves with the same cohomology classes, which are graphically represented by dotted lines. There are $h^{4}(X)$ such graphs of type $C$, $D$ and $E$. 
 
Adding up all these contributions from the graphs above, we find a lower bound for the dimension of $H^{2}(\fscheme)$. As promised, this coincides with the number given by equation $\eqref {upper}$. This completes the proof of theorem $\ref {diva}$. 

It remains to explain why the five types of fixed loci listed above have exactly one negative weight on their normal bundles. Needless to say, the computation involves the Konsevich-Graber-Pandharipande recipe for computing the weights. This is a straightforward argument for the most part. There are two ingredients which are important to the count of the negative weights. 

\begin {enumerate}
\item [(a)] For each curve $R=R_i$, the tangent space $T_{P} R$ at the fixed point $P$ has one positive weight, as all weights on $T_P X$ are positive. Therefore the tangent space $T_{q_{i}}R$ has the opposite/negative weight. When $R$ is the curve joining $q_i$ to $q_{ij}$, the tangent space at $q_i$ has one positive weight: there should be only one negative weight on $T_{q_i}X$ which we've seen occurs along the curve $R_i$. Therefore the tangent space of $R$ at $q_{ij}$ has the opposite/negative weight. Moreover, an explicit computation of the weights shows that each vertex of valency $2$ in the graphs of type $D, E$ contributes with the positive weight $\omega_{\mathfrak f_1}+\omega_{\mathfrak f_2}$. 
\item [(b)] For each rational curve $R=R_{j}$ joining $P$ to $q_j$ and each map $f: \p^1\to R$ of degree $d$ which is totally ramified over $P$ and $q_{j}$, the number of negative weights on $H^{0}(\p^1, f^{\star}TX)$ equals $d$.  The case in hand can be checked rather easily recalling the following "Euler" sequence on $X$: 
\begin {equation}\label {eulerflag} 0 \to TX \to \bigoplus_{i} Hom(\mathcal S_i, \mathcal Q_{i})\to \bigoplus_{i} Hom (\mathcal S_i, \mathcal Q_{i+1})\to 0. 
\end {equation}
It is enough to compute the weights on the virtual $\cs$ representation $$\bigoplus_{i} H^0(f^{\star}(\mathcal S_i^{\star}\otimes \mathcal Q_{i}))-\bigoplus_{i}H^{0}(f^{\star}(\mathcal S_i^{\star}\otimes \mathcal Q_{i+1})).$$ We claim that when $i\neq j$, there are no negative weights on $H^{0}(f^{\star}(\s^{\star}_i \otimes \q_i))$. Indeed, we observe that for $i\neq j$, we have the equivariant isomorphism $$f^{\star} \s_i= W(m_i) \otimes \o$$ so the weights on $H^{0}(f^{\star}(\s^{\star}_i \otimes \q_i))$ are the positive numbers $\lambda_{u}-\lambda_{v}$ for $v \leq m_i < u$. 

We claim $d$ negative weights on $H^{0}(f^{\star}(\s^{\star}_j \otimes \q_j))$. Equation $\eqref {line}$ shows that we have an equivariant identification: \begin {equation} \label {neg2} \s_j|_{R_j} = W(m_j-1) \otimes \o \oplus \o(-1).\end {equation} In the above, the $\cs$ action on $\o(-1)$ has the weight $\lambda_{m_j}$ at $[1:0]$ and the weight $\lambda_{{m_j}+1}$ at $[0:1]$.  

The exact sequence $\eqref {tautflag}$ tensored with $f^{\star} \s_j^{\star}$ shows that $H^{0}(f^{\star} \s^{\star}_j \otimes f^{\star} \q_j)$ is equivalent to the $\cs$ virtual representation \begin {equation} \label {neg1} V \otimes H^{0} (f^{\star} \s^{\star}_j) - H^{0} (f^{\star} \s^{\star}_j \otimes f^{\star} \s_j) + H^{1} (f^{\star} \s^{\star}_j \otimes f^{\star} \s_j).\end {equation} We can now consider each of the three terms in $\eqref {neg1}$ individually. Using $\eqref {neg2}$, we rewrite the first term as $$V \otimes W(m_j-1)^{\star} \oplus V \otimes H^{0} (\mathcal O_{\p^1}(d)).$$ The first summand has weights $\lambda_u - \lambda_v$ for $v\leq m_j-1$ and all $u$. In addition, we also have the weights $\lambda_u - \frac{1}{d}(a\lambda_{m_j}+b\lambda_{m_j+1})$ for all nonnegative $a, b$ such that $a+b=d$. The second term in $\eqref {neg1}$ can be rewritten as $$W(m_j-1)^{\star} \otimes W(m_j-1) \oplus H^{0}(\o(d)) \otimes W(m_j-1) \oplus \c$$ with a trivial action on the last term. This has exactly the same non-zero weights as the first term in $\eqref{neg1}$, except that we need to require that $u\leq m_j-1$. Finally, the third term in $\eqref {neg1}$ equals $$W(m_j-1)^{\star}\otimes H^{1}(\o(-d))$$ with positive weights $-\lambda_{u}+ \frac{a\lambda_{m_j}+b\lambda_{m_j+1}}{d}$ for all $u\leq m_j-1$ and $a,b$ are positive integers summing up to $d$. Summarizing, we find that the only negative weights in the list above are the $d$ values $\frac{a}{d} (\lambda_{m_j}- \lambda_{m_j+1})$ for $1\leq a\leq d$. 

In the same way we verify that there are no negative weights on $H^{0}(f^{\star}(\s_i^{\star} \otimes \q_{i+1}))$. We leave the details to the reader since they are similar to the computations above. Our initial claim is now proved. 

Similarly, one shows that the number of negative weights on $H^{0}(R, TX)$ is exactly $2$ whenever $R$ is a rational curve in the cohomology class dual to $\beta_j$ of one of the following types: \begin {itemize} \item the rational curve joining $P$ to one of the points $r_j$ or $r'_j$. \item the rational curve  joining $q_i$ to $q_{ij}$. \item the rational curve joining $q_{j+1}$ to $s_j$ or the curve joining $q_{j-1}$ to $s'_j$. \end {itemize}

\end {enumerate}

\begin {remark} The author believes that the arguments of the last few subsections can be repeated for $X=G/P$. The proofs should not be any more difficult. The count of negative $\cs$ weights on various cohomology groups can be done in terms of the roots of $\mathfrak p$. Moreover, the count of the negative weights on the various tangent spaces follows from standard computations. We leave these arguments to the interested reader.\end {remark}

\section {The codimension 2 classes.} This section contains the proof of proposition $\ref {eval}$ which will be instrumental in the final section of this paper. In this section, {\bf we will only consider the interesting case $d>1$.}

It is useful to introduce the following piece of {\bf notation} for a collection of classes which will appear frequently in our computations. First, we make the convention that $a,b$ will always denote indices adding up to $d$. Next, we consider $S_1, S_2$ two disjoint, possibly empty subsets of $\{1, 2, \ldots, n\}$ and we let $C_1, C_2$ be two possibly empty collections of cohomology classes on $X$. We will denote by $\Delta_{a, b}(S_1, S_2 | C_1, C_2)$ the class on $\fscheme$ represented by the subscheme of stable maps with two components of degrees $a$ and $b$ carrying the markings $S_1$ and $S_2$ respectively. Moreover, we require that these components intersect generic subvarieties in the cohomology classes $C_1$ and $C_2$ respectively. The sum of all cohomology classes $\Delta_{a,b}(S_1, S_2|C_1, C_2)$ for all possible degrees is denoted by $\Delta(S_1, S_2 | C_1, C_2)$. We denote by $\bullet$ the empty sets among $S_1, S_2, C_1, C_2$. If $C_1, C_2$ are both empty, we will omit writing both bullets. 

\begin{proposition}\label{eval} For all $1\leq i\neq j\leq n$ and for all codimension $2$ classes $\alpha$ on the Grassmannian $\bg$, the following relation is true:
\begin {equation}\label {re2}
ev_i^{\star}\alpha-ev_j^{\star}\alpha-\psi_j \kappa(\alpha)=- \Delta(\{i\}, \{j\} | \alpha, \bullet). 
\end {equation}
\end {proposition}

As a first reduction, we observe that the equations $\eqref {re2}$ pull back compatibly under the forgetful maps $\pi:\overline M_{0,n+1}(\bg, d)\to \overline M_{0,n}(\bg, d)$. This is an immediate consequence of the relation $\pi^{\star}\psi_j=\psi_j-\Delta_{0, d}(\{j,n+1\},\{1, \ldots, j-1, j+1, \ldots, n\})$. Henceforth, we will assume $n=2$. Secondly, we observe that it is enough to work on Grassmannians of rank $r\geq 3$ subspaces since the rank $2$ case will follow by restriction to the space of maps to ${\bf G}(2, V)\hookrightarrow \bg$.

As mentioned in the introduction, the strategy of the proof will be to enumerate all codimension $2$ cycles and compare their number with the dimension of the Chow group which will be computed by a combination of localization and Deligne spectral sequence techniques. As a warm up exercise, we will undertake this task in the next subsection for the space of maps to $\pr$. 

\subsection {Maps to Projective Spaces.} To begin with, we consider the case of maps without markings $\overline M_{0,0}(\pr, d)$. There are several cases to consider: $r=1$, $r=2$ and $r\geq 3$.

We will consider $X=\p^1$ first. In \cite {O1} and \cite {O2}, we saw that the whole cohomology of the moduli space of stable maps to $\p^1$ is generated by the tautological classes. The search for the codimension $2$ tautological classes gives the following results: \begin {enumerate} \item boundary classes of nodal maps with three irreducible components, \item classes of nodal stable maps whose node maps to a fixed point in $\p^1$. \end {enumerate}

It turns out the classes in (1) and (2) are in fact linearly independent. To see this it is enough to match their number with the dimension of $H^{4}(\overline M_{0,0}(\p^1, d))$. This dimension count can be done in several ways - either using Deligne's spectral sequence or localization. 

Let us briefly indicate how the localization argument works, then turning to the computation by Deligne's spectral sequence. We will use a $\cs$ action on $\p^1$ with two fixed points which we call $0, 1$ such that the tangent bundle at $0$ has a positive weight, while the tangent bundle at $1$ has a negative weight.  As before, we index the fixed point loci on $\overline M_{0,0}(\p^1, d)$ by decorated graphs $\Gamma$. 

In \cite {O2}, we counted the negative weighs on the normal bundles of the fixed loci (see lemma $3$ in that paper). More precisely, we let $\mathfrak s$ denote the number of vertices of $\Gamma$ labeled by $1$ with at least three incident flags and let $\mathfrak u$ denote the number of vertices of $\Gamma$ labeled by $1$ with one incident flag. The number of negative weights on the normal bundle of the fixed locus indexed by $\Gamma$ is $d - \mathfrak u + \mathfrak s$. An easy argument shows that there are $4$ types of graphs indexing fixed loci with at most $2$ negative weights on their normal bundle. These graphs are shown in the figure below. The number of negative weights is indicated in a box to the left of each graph and the degrees of the edges are written below the graphs.

\begin {figure}[htb]\label{2p1}
\scalebox {.8}{\psfig{figure=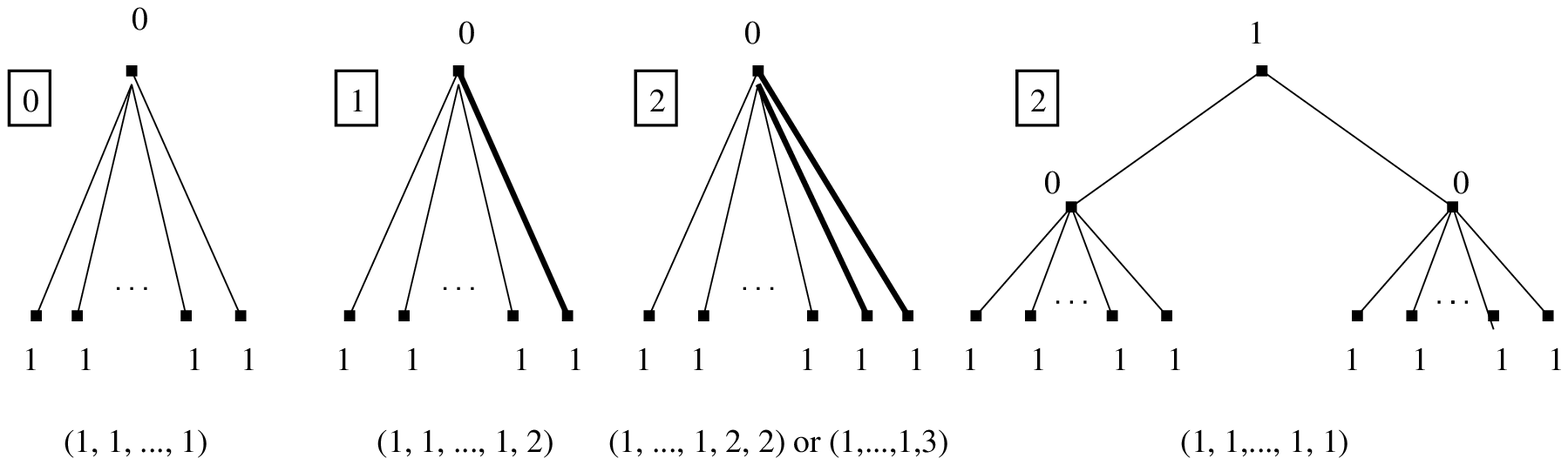}}
\caption {The fixed loci for the $\cs$ action on $\overline M_{0,0}(\p^1, d)$.}
\end {figure}

The first graph gives the fixed locus $\overline M_{0,d}/S_d$ with no negative weights on the normal bundle, while the second graph gives the fixed locus $\overline M_{0, d-1}/S_{d-2}$ with one negative weight. There are $2+[d/2]$ graphs with $2$ negative weights. We find that the dimension of $H^4$ equals $$h^{4}({\overline M}_{0,d})^{S_d} + h^{2} (\overline M_{0, d-1})^{S_{d-2}}+ 2 + [d/2] = h^{4}({\overline M}_{0,d})^{S_d} + d-2 + [d/2]$$ where we also used lemma $\ref {dim}$ to evaluate the second term of the sum.

It remains to compute the first term. By a well know theorem of Keel, we already know that the generators of $H^{4}(\overline M_{0,d})^{S_d}$ are the boundary classes of curves with at least $3$ irreducible components. These boundary classes $\mathcal {B}_{ijl}$ are indexed by triples $(i, j, l)$ such that $i+j+l = d, \; j\geq 1,\; 2\leq i\leq l,$ these integers corresponding to the number of marked points on each component. 

{\bf Claim.}
The classes $\mathcal B_{ijl}$ form a basis of $H^{4}(\overline M_{0,d})^{S_d}$. \vskip.06in

We show these classes are linearly independent by matching their number with the actual dimension of $H^4(\overline M_{0,d})^{S_d}$. The dimension computation is identical to that in lemma $\ref {dim}$ by averaging out the traces of all $\sigma\in S_d$ on $H^{4}(\overline M_{0,d})$. We will omit the details. An alternate proof will be established by the arguments below. 

We will now redo the same computation making use of Deligne's spectral sequence. There are two cases to consider depending on the parity of $d$. Let us show the details for $d=2k$. An identical argument also works for odd $d$'s - the numerical details are slightly different. We make use of the fact that the boundary divisors in the space of stable maps have normal crossings. However, these boundary divisors have self-intersections and writing down the Deligne spectral sequence with the right system of coefficients is delicate, as we have to account for automorphisms and self-intersections; the relevant details were explained in the first section of \cite {O1}.

In our case, the $k$ boundary divisors $D_1, \ldots D_k$ correspond to nodal maps whose degrees on the components are $(1, 2k-1)$, \ldots $(k,k)$. The codimension two strata are denoted by $D_{ijl}$ for $1\leq i\leq l$ and $1\leq j$ with $i+j+l=d$; they correspond to stable maps with three components, the degree on the middle component being $j$. For $1\leq i < k$, the stratum $D_{i, 2k-2i, i}$ is the self intersection of $D_{i, 2k-i}$; there are no anti invariant classes in its zeroth cohomology group because of the $\mathbb Z/2\mathbb Z$ symmetry which switches the edges. Therefore, this term does not appear in the Deligne spectral sequence. We obtain the following complex: 
\begin {equation}\label {del}\bigoplus_{1\leq i<l, 1\leq j, i+j+l=d} H^{0}(D_{ijl}) \to \bigoplus_{i=1}^{k-1} H^{2}(D_i)\oplus H^{2}(D_k)^{-} \to H^{4}(\overline M_{0,0}(\p^1,d)).
\end {equation}
Here the minus superscript on $H^{2}(D_k)^{-}$ stands for the anti invariant part of the cohomology under the sign representation of $Aut(\Gamma)$, as explained in \cite {O1}. Recall that each automorphism has a sign, given by the action it induces on the determinant $\det (\text { Edge }(\Gamma))$. The minus sign refers to cohomology classes which are invariant under this sign representation. Since the edge of the dual graph indexing $D_k$ is preserved under the $\mathbb Z/{2\mathbb Z}$ symmetry of the graph, $H^{2}(D_k)^{-}$ is the same as the $\mathbb Z/{2\mathbb Z}$ invariant part of $H^{2}(\overline M_{0,1}(\p^1, k)\times_{\p^1} \overline M_{0,1}(\p^1, k))$. We use the computation of lemma $\ref {dim}$ to conclude that the dimension of the middle term of the complex $\eqref {del}$ is $2k^{2}-2k+1$. The first term is easily seen to be $k^{2}-2k+1$ dimensional. 

We claim the dimension of $H^{4}$ is $k^2$, which also turns out to be the number of generators we have exhibited for $H^{4}$. It suffices to show that the alternating sum of dimensions of the terms in the complex $\eqref{del}$ is $0$. This alternating sum equals the coefficient of $q^{2\text {dim } - 4}$ in the virtual Poincare polynomial of $M=M_{0,0}(\p^1,d)$. By definition, the Poincare polynomial can be computed from the associated graded of the Hodge weight filtration: \begin {equation} \label {hodge} P(M)=\sum_{i,j} (-1)^{i+j} Gr_j^{W}(H^{i}_c(M)) q^{j}\end {equation} In \cite {M1}, it is proved that $P(M)=q^{4d-4}$. The claim follows from these observations.

The argument above is no more complicated when we deal with general projective spaces. We will pick a $\cs$ action with weights $\lambda_0,\ldots, \ldots, \lambda_r$ such that $\lambda_i$ is much bigger than $\lambda_{i-1}$. We denote by $0, 1, \ldots, r$ the isolated fixed points of this action. To find the dimension of $H^4$ we use localization. In addition to the fixed loci listed above for $\p^1$ we have $4$ more types of graphs when $r=2$ and one additional graph for $r=3$. There are no new graphs added for $r>3$. This is an aspect of the "stabilization of cohomology" theorem proved in \cite {BH}. Our computation shows that when $r=2$ we gain $d$ dimensions (we will need to invoke lemma $\ref {dim}$ again to compute the contribution of the first graph). For $r\geq 3$ we gain $d+1$ dimensions. 

\begin {figure}[htb]
\scalebox {.8}{\psfig{figure=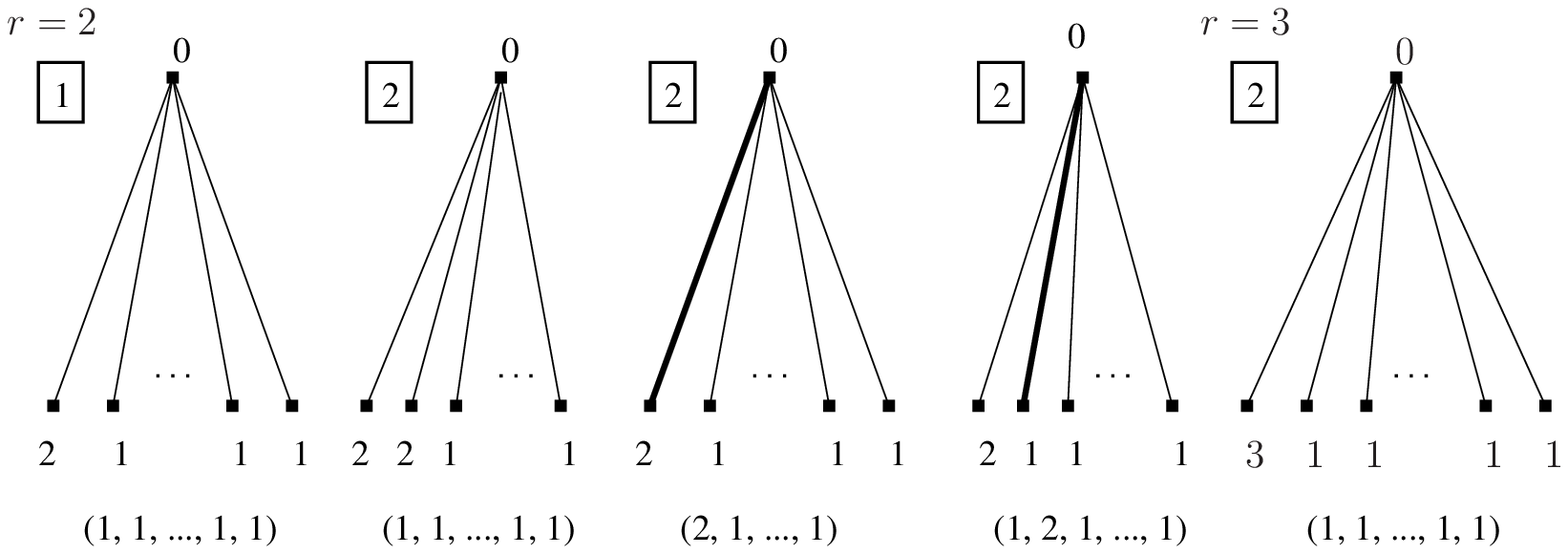}}\label{2pr}
\caption {The fixed loci for the $\cs$ action on $\overline M_{0,0}(\pr, d)$.} 
\end {figure}

We find generators for $H^4$ by enumerating all tautological classes and matching their number with the dimension we just computed. We obtain the following:

\begin {proposition}\label {codimension} {\bf A.}
The following codimension $2$ classes on $\overline M_{0,0}(\p^r, d)$ form a basis for the codimension $2$ Chow group. 
\begin {itemize} 
\item [(A.1)] the boundary classes of maps whose domain has at least three components,
\item [(A.2.1)] the nodal classes of maps whose node is mapped to a codimension $1$ subspace, 
\item [(A.2.2)](when $r\geq 2$) classes of nodal maps, one component passing through a fixed codimension $2$ subspace, 
\item [(A.3.1)](when $r\geq 2$) classes of maps whose images pass through two general codimension two subspaces,
\item [(A.3.2)](when $r\geq 3$) the class of maps intersecting a codimension $3$ subspace. 
\end {itemize}
\end {proposition}

We will consider now the case of one marked point. Enumeration of the codimension $2$ tautological classes yields the following results:
\begin {enumerate}
\item [(B.1)] classes of stable maps with at least three irreducible components. 
\item [(B.2.1)] classes of stable maps with two components such that the marked point maps to a fixed hyperplane in $\pr$.
\item [(B.2.2)] classes of stable maps with two components, one component intersecting a fixed codimension $2$ subspace in $\pr$ (for $r\geq 2$).
\item [(B.3.1)] classes of maps whose marked point maps to a fixed codimension $2$ subspace of $\pr$ (if $r\geq 2$). 
\item [(B.3.2)] classes of maps whose marked point maps to a fixed hyperplane and whose images intersect a codimension $2$ subspace (if $r\geq 2$).
\item [(B.3.3)] classes of maps which intersect two general codimension $2$ subspaces of $\pr$ (if $r\geq 2$).
\item [(B.3.4)] classes of maps which intersect a codimension $3$ subspace of $\pr$ (if $r\geq 3$).
\end {enumerate}
The reader may easily observe that we omitted: \begin {itemize}\item the classes of stable maps with two irreducible components, whose node maps to a fixed hyperplane in $\pr$. \end {itemize} By equation $\eqref {evsum}$ below, these can be expressed in terms of classes in $(B.1)$, $(B.2.1)$ and $(B.2.2)$. 

We invite the reader to compute the actual dimension of the cohomology using Deligne's spectral sequence and observe that in fact we need to have one relation between the generators above when $r\geq 2$, whereas for $r=1$, the {\it nonzero} classes above are in fact independent. The missing relation for $r\geq 2$ was explained in lemma 5 of \cite {O1}. We proved there that \begin {equation}\label {marked} ev_1^{\star} H^2 - \frac{1}{d}ev_1^{\star} H\cdot  \kappa(H^2)+\frac{1}{d^2} \kappa(H^2)^2 - \frac{1}{d} \kappa(H^3)\end {equation} is a sum of classes of type $(B.1)$ and $(B.2)$. 

The case of two markings is entirely similar and we leave the details to the interested reader. We only mention that in this case we need to account for more than one relation. However, all these relations can be explained by lemmas 1.1.1, 1.1.2 in \cite {divisors}, equations $\eqref {diff}$ - $\eqref {evsum}$ below and the relation $\eqref {marked}$ above. We summarize our findings:

\addtocounter{proposition}{-1}
\begin {proposition}
{\bf B.} The classes $(B.1)$, $(B.2.1)$, $(B.2.2)$, $(B.3)$ generate the codimension $2$ Chow group of $\schemea$. The only possible relation between the {\bf nonzero} classes above is $\eqref {marked}$ (for $r\geq 2$).\\
\end {proposition}
\addtocounter{proposition}{-1}
\begin {proposition}
{\bf C.} The following classes, when nonzero, form a basis for the codimension $2$ Chow group of $\schemeb$: 
\begin {enumerate}
\item [(C.1)] classes of maps with at least three components.
\item [(C.2.1)] classes of maps with two components, the markings $1$ and $2$ being on different components, one of the components intersecting a codimension $2$ subspace in $\pr$ (for $r\geq 2$).  
\item [(C.2.2)] classes of maps with two components, the markings $1$ and $2$ being on different components, the node mapping to a fixed hyperplane.
\item [(C.2.3)] the class of maps with two components, one of degree $0$ containing both markings and mapping to a fixed hyperplane. 
\item [(C.2.4)] the class of maps with two components, one of these components containing both markings, the other intersecting a fixed codimension $2$ subspace of $\pr$ (if $r\geq 2$).
\item [(C.3)] classes of maps such that one marking maps to a codimension $2$ subspace of $\pr$ (if $r\geq 2$).
\end {enumerate}
\end {proposition}

\subsection {A useful relation.} In this subsection we will prove a particular case of $\eqref {re2}$ which will be important in deriving the general result.  
 
\begin {lemma}\label{re2p}
Proposition $\ref {eval}$ holds on the space of maps to $\pr$. That is, the following relation is satisfied: \begin {equation}\label{2m}ev_1^{\star}H^2-ev_2^{\star}H^2 - \psi_2 \kappa(H^2)=-\Delta(\{1\}, \{2\} | H^2, \bullet).\end {equation}
\end {lemma}

{\bf Proof.} The proof of this result is a consequence of the divisorial relations below. These are theorem $1$ in \cite {LP}, and lemma 2.2.2 in \cite {divisors} respectively. 

\begin {equation}\label {diff}
ev_1^{\star} H - ev_2^{\star} H = d\psi_2 - \sum a \Delta_{a,b} (\{1\}, \{2\}).
\end {equation}
\begin {equation} \label {psisum}
\psi_1 + \psi_2 = \Delta(\one, \two)
\end {equation}
\begin {equation} \label {strange}
\psi_i+\frac{2}{d} ev_i^{\star} H - \frac{1}{d^2} \kappa(H^2) = \sum \frac{b^2}{d^2}\Delta_{a,b}(\{i\}, \b) 
\end {equation}
As a consequence of $\eqref {psisum}$ and $\eqref {strange}$ we obtain the following two equations: \begin {equation} \label {psi2} \psi_2=\frac{2}{d} ev_1^{\star} H - \frac{1}{d^2}\kappa(H^2) + \text {boundaries } \end {equation}
\begin {equation}\label {evsum} ev_1^{\star} H + ev_2^{\star} H = \frac{1}{d} \kappa(H^2)  -\frac{ab}{d} \Delta_{a,b}(\one, \two) + \frac{b^2}{d} \Delta_{a,b}(\{1,2\}, \b)\end {equation} 
We multiply $\eqref{diff}$ by $\eqref {evsum}$ to obtain the following expression for $ev_1^{\star}H^2-ev_2^{\star} H^2$: \begin {equation}\label {squares} \left(d\psi_2-\sum a\Delta_{a,b}(\one, \two)\right)\left(\frac{1}{d} \kappa(H^2) - \frac{ab}{d} \Delta_{a,b}(\one, \two) + \frac{b^2}{d}\Delta_{a,b}(\{1,2\}, \b)\right)\end {equation} This product has the form $\psi_2 \cdot \kappa(H^2)+\text {boundary terms }$. We observe that $ev_1^{\star}H^2-ev_2^{\star} H^2-\psi_2 \cdot \kappa(H^2)$ restricts to $0$ on the moduli space of maps to any $\p^1\hookrightarrow\pr$. This observation and the independence result in proposition $\ref {codimension}$.${\bf C}$ allows us to conclude that among the boundary terms in the expansion of $\eqref {squares}$ we cannot have classes of type $(C.1)$, $(C.2.2)$, $(C.2.3)$. We will therefore work modulo such boundary terms when evaluating the product $\eqref {squares}$. 

Such codimension $2$ boundary terms are obtained for example by multiplying out two distinct divisors $\Delta_{a,b}$. Thus, we will only need to worry about self intersections of such divisor classes. Similarly, the product $\psi_2 \cdot \Delta_{a,b}(\{1,2\}, \b)$ can be rewritten as sum of classes of type $(C.1)$ since the $\psi$ class on the component with three special points is a boundary class (\cite {divisors}, lemma 1.1.1).  We conclude that: \begin {eqnarray*} ev_1^{\star}H^2&-&ev_2^{\star} H^2- \psi_2 \cdot \kappa(H^2)= -ab \psi_2 \cdot \Delta_{a,b}(\one, \two) + \frac{a^2 b}{d}\Delta_{a,b}(\one, \two)^2 -\\&-&\frac{a}{d}\Delta_{a,b}(\one, \two|\b, H^2) - \frac{a}{d}\Delta_{a,b}(\one, \two|H^2, \b).\end {eqnarray*} We will write $\mathcal D_{a,b}$ for the class of nodal maps of degrees $a$ and $b$, with markings $1$ and $2$ on the first and second component respectively, such that the node maps to a fixed hyperplane of $\pr$. Using $\eqref {psi2}$ we obtain that $$\psi_2 \cdot \Delta_{a,b}(\one, \two)=\frac{2}{b}\mathcal D_{a,b} - \frac{1}{b^2} \Delta_{a,b}(\one, \two|\b, H^2)+\text {boundaries }.$$ The divisor $\Delta_{a,b}(\one, \two)$ is dominated by $\overline M_{0,\one\cup\{\star\}}(\pr, a)\times_{\pr}\overline M_{0, \{\bullet\}\cup \two}(\pr, b)$ and its normal bundle has Chern class $-\psi_{\star}-\psi_{\b}$. Using the self-intersection formula and $\eqref {strange}$ we obtain that modulo boundaries the following is true: $$\Delta_{a,b}(\one, \two)^2= \frac{2}{a}\mathcal D_{a,b} - \frac{1}{a^2} \Delta_{a,b}(\one, \two|H^2, \b) + \frac{2}{b} \mathcal D_{ab} - \frac{1}{b^2} \Delta_{a,b}(\one, \two|\b, H^2).$$ The lemma follows easily from the last three equations. 

\begin {remark}\label {red}As a consequence, we prove proposition $\ref {eval}$ for $\alpha = c_1(\det \q)^2$ by making use of the map $\overline M_{0,2}(\bg, d)\hookrightarrow \overline M_{0,2}(\pr, d)$ induced by the Plucker embedding. \end {remark}

\subsection {Maps to Grassmannians.} We will repeat the reasoning of subsection $3.1$ in the context of the unmarked stable map spaces to Grassmannians. However, there is no result similar to Manin's computation of the Hodge polynomial of the open part. Instead, we will appeal to localization. There are complications arising in this approach since we do not have a formula counting the number of negative weights on the normal bundles of the fixed loci, as we did for $\pr$. Instead, we will use localization to obtain a lower bound for the dimension of the codimension $2$ Chow group of $\overline M_{0,0}(\bg, d)$. This bound will turn out to be sharp, as an easy enumeration of the codimension $2$ tautological classes will show. We will combine this with a formal computation involving Deligne's spectral sequence to find the dimension of the cohomology when we increase the number of marked points. By comparing this with the number of generators, we will be able to derive the existence of one relation, which will be explicitly identified next.

We assume $\bg$ is the Grassmannian of quotients of dimension $r$ of some $N$ dimensional vector space $V$ such that $r\geq 3$ and $N - r\geq 3$. We will use a $\cs$ action on $V$ with weights $\lambda_1 < \ldots< \lambda_N$ as in section $2.4$. We label the fixed points of the induced $\cs$ action on $\bg$ by the number of positive weights on their tangent bundle. There will be one fixed point labeled $0$, one fixed point labeled $1$, two fixed points labeled $2$ and (since $r\geq 3$) there will be $3$ fixed points labeled $3$. Moreover, all these fixed points can be joined by a degree $1$ rational curve to the fixed point $0$. 

We will consider the induced $\cs$ action on $\overline M_{0,0}(\bg, d)$ and exhibit fixed loci with at most two negative weights on their normal bundles. In fact, we have already done most of the work for projective spaces. We already know {\it all} such graphs for $\p^3$ as a result of the previous subsections. These graphs are exhibited in figures $3$ and $4$. They still make sense for $\bg$, but their number will be higher in this case, since there are multiple fixed points on $\bg$ with the same labels. A moment's though shows that the five graphs we exhibited in figure $4$ contribute at least $h^{2}(\overline M_{0,d})^{S_{d-1}}+ 2 + 1 + 1+ 2 = d+3$ extra dimensions to $h^{4}(\overline M_{0,0}(\bg, d))$ over $h^{4}(\overline M_{0,0}(\p^3, d))$. Presumably, there could be more graphs which contribute, so we can only conclude: \begin {equation}\label {comp}h^{4}(\overline M_{0,0}(\bg, d))-h^{4}(\overline M_{0,0}(\p^3, d))\geq d+3.\end {equation}

We will now argue that the opposite inequality is true. We will invoke the main result of \cite {O1} which shows all cohomology of $\overline M_{0,0}(\bg, d)$ is tautological. Enumeration of the codimension $2$ cohomology classes yields:
\begin {itemize}
\item [(D.1)] the classes we exhibited for $\p^3$ in proposition $\ref {codimension}$.{\bf A}, but replacing the arbitrary hyperplane in $\p^3$, the codimension $2$ subspace and the codimension $3$ subspace used there by smooth subvarieties in the cohomology classes $c_1(\q)$, $c_1(\q)^2$ and $c_1(\q)^3$ respectively.
\item [(D.2)] boundary classes of nodal maps, one of the components intersecting a fixed subvariety in the class $c_2(\q)$. 
\item [(D.3.1)] the class of maps which intersect two generic subvarieties both in the cohomology class $c_2(\q)$.
\item [(D.3.2)] the class of maps which intersect two generic subvarieties in the cohomology classes $c_1^2(\q)$ and $c_2(\q)$.
\item [(D.3.3)] the class of maps which intersect a subvariety in the cohomology class $c_1(\q) c_2(\q)$. 
\item [(D.3.4)] the class of maps which intersect a subvariety in the cohomology class $c_3(\q)$. 
\end {itemize} 

This enumeration is enough to establish that $\eqref {comp}$ is an equality. We find the following result. The independence of the degree seems to be true for all flag varieties; compare to the similar statement for projective spaces proved in \cite {O1} and to the similar remark in \cite {M1}. \addtocounter{proposition}{-1}\begin{proposition} {\bf D.} The classes listed above, if nonzero, form a basis for the complex codimension $2$ cohomology of $\overline M_{0,0}(\bg, d)$. The lowest piece of the Hodge structure on the open part is generated by the classes $(D.3)$ and $(A.3)$. Its dimension is degree independent.\end {proposition} We now compare the virtual Poincare polynomials $P(M_{0,0}(X, d))$ defined in $\eqref {hodge}$ when $X$ is either $\p^3$ or $\bg$. We use Deligne's spectral sequence. There are two cases two consider depending on the parity of $d$. Let us only spell out the case when $d=2k$. We borrow the notation of subsection $3.1$ and we write $D_{ijl}^{X}$ for the class of maps to $X$ with three components of degrees $i, j, l$ and $D_i^X$ for the boundary class of maps with two components of degree $i$ and $2k-i$ on the space $\overline M_{0,0}(X, d)$. We write $D(X)$ for the complex dimension of $\overline M_{0,0}(X, d)$. Just as in $\eqref {del}$, we obtain complexes 
\begin {equation}\label {ho1}
0\to \bigoplus_{i=1}^{k}H^0(D_i^X)\to H^{2}(\overline M_{0,0}(X, d))\to 0
\end {equation}

\begin {equation}\label {ho2}
0\to \bigoplus_{1\leq i<l, j\geq 1, i+j+l=d} H^{0}(D_{ijl}^X)\to \bigoplus_{i=1}^{k-1} H^{2}(D_i^X)\oplus H^{2}(D_k^X)^{-}\to H^{4}(\overline M_{0,0}(X, d))\to 0.
\end {equation}
The alternating sums of the dimensions of the terms in $\eqref {ho1}$ is the coefficient $a_2(X)$ of $q^{2D(X)-2}$ in the virtual Poincare polynomial $P(M_{0,0}(X,d))$, while for the complex $\eqref {ho2}$ the alternating sums of dimensions is the coefficient $a_4(X)$ of $q^{2D(X)-4}$ in the same polynomial. Comparing the complexes $\eqref {ho1}$ for $X=\p^3$ and for $X=\bg$ using the dimension formula $\eqref {upper}$ we obtain \begin {equation}\label {2ndcoeff} a_2(\bg)=a_2(\p^3)+1.\end {equation} Now, we can use $\eqref {comp}$ to compare the last term of the complex $\eqref {ho2}$. The dimension formula $\eqref {upper}$ shows that the dimensions of the middle terms differ by $2$ when $1\leq i\leq k-1$ and by $1$ for $i=k$. We conclude that: \begin {equation}\label {4thcoeff}a_4(\bg)=a_4(\p^3)+(d+3)-(2(k-1)+1)=a_4(\p^3)+4.\end {equation} The map $\pi: M_{0,1}(X, d)\to M_{0,0}(X, d)$ is a locally trivial fibration with fibers isomorphic to $\p^1$, therefore: $$P(M_{0,1}(X, d))=(q^2+1) P(M_{0,0}(X, d)).$$
Writing ${\overline a}_4(X)$ for the coefficient of $q^{2(D(X)+1)-4}$ in $P(M_{0,1}(X, d))$, and using the above equation together with $\eqref {2ndcoeff}$, $\eqref {4thcoeff}$ we obtain that: \begin {equation}\label {coeff}{\overline a}_4(\bg)={\overline a}_4(\p^3)+5.\end {equation} The same conclusion is valid when $d$ is odd.

We now repeat this argument "backwards" for the moduli spaces $\overline M_{0,1}(X, d)$. Using the complex $\eqref {ho2}$, equation $\eqref {upper}$ to compare the middle terms, and equation $\eqref {coeff}$ to compare the alternating sums of dimensions we find: $$h^{4}(\overline M_{0,1}(\bg,d))=h^{4}(\overline M_{0,1}(\p^3, d))+ 2d+3.$$ We again list the possible tautological generators for the codimension $2$ Chow group:
\begin {itemize}
\item [(E.1)] All classes we exhibited for $\p^3$ in proposition $\ref {codimension}$.{\bf B}, but replacing the arbitrary hyperplane in $\p^3$, the codimension $2$ subspace and the codimension $3$ subspace there by smooth subvarieties in the cohomology classes $c_1(\q)$, $c_1(\q)^2$ and $c_1(\q)^3$.
\item [(E.2)] boundary classes of nodal maps, one of the components intersecting a fixed subvariety in the class $c_2(\q)$. 
\item [(E.3.1)] the class of maps whose marked point maps to a subvariety in the cohomology class $c_2(\q)$.
\item [(E.3.2)] the class of maps whose marked point maps to a subvariety in the cohomology class $c_1(\q)$ and whose images intersect a subvariety in the cohomology class $c_2(\q)$.
\item [(E.3.3)] the class of maps which intersect two generic subvarieties both in the cohomology class $c_2(\q)$.
\item [(E.3.4)] the class of maps which intersect two generic subvarieties in the cohomology classes $c_1^2(\q)$ and $c_2(\q)$.
\item [(E.3.5)] the class of maps which intersect a subvariety in the cohomology class $c_1(\q) c_2(\q)$. 
\item [(E.3.6)] the class of maps which intersect a subvariety in the cohomology class $c_3(\q)$. 
\end {itemize} 

We now observe that the span of the classes in $(E.1)$ is exactly $h^{4}(\overline M_{0,1}(\p^3, d)))$ dimensional. Indeed, we only need to account for a relation analogous to $\eqref {marked}$, replacing $H$ by $c_1(\q)$. Such a relation is obtained from the one for $\pr$ by making use of the Plucker embedding. Restricting to a generic $\p^3\hookrightarrow\bg$ and using proposition $\ref {codimension}$.{\bf B.} we see that there cannot be any more relations. 

This observation and formula $\eqref {coeff}$ shows that there should be {\it exactly} one more relation involving classes in $(E.2)$ and $(E.3)$. Write $c_i$ for the Chern classes $c_i(\q)$. For constants $\gamma_i^{(d)}$, not all zero, we conclude that \begin {equation}\label{1m} \gamma_1^{(d)} \cdot ev_1^{\star} c_2+ \gamma_2 ^{(d)}\cdot ev_1^{\star} c_1 \cdot \kappa(c_2)+\gamma_3^{(d)} \cdot \kappa (c_2)^2 + \gamma_4^{(d)} \cdot \kappa(c_1^2, c_2) +\gamma_5^{(d)} \cdot \kappa(c_1c_2)+ \gamma_6^{(d)} \cdot \kappa(c_3) +\end {equation} $$+\gamma_7^{(d)} \cdot ev_1^{\star}c_1^2+\gamma_8^{(d)} \cdot ev_{1}^{\star} c_1 \cdot\kappa(c_1^2)+\gamma_9^{(d)}\cdot \kappa(c_1^2)^2+\gamma_{10}^{(d)} \cdot \kappa(c_1^3)$$ is sum of boundary classes in $(E.2)$ and {\it boundary} classes in $(E.1)$. 

Next, we will identify the constants $\gamma_i^{(d)}$. The first observation is that the $\gamma_i^{(d)}$'s do not depend on the dimension of $V$. Indeed, assume that for two vector spaces $V\hookrightarrow W$ we had different coefficients for the relations $\eqref {1m}$. Restricting the equation for ${\bf G}(r, W)$ to ${\bf G} (r, V)$ we obtain two distinct relations, which is impossible. Therefore, to find the constants $\gamma_i^{(d)}$ we can assume that the dimension of $V$ is as large as we want. 

Our strategy will be to pull back equation $\eqref {1m}$ to moduli spaces of maps to $\p^3$. The vector bundle $$\mathcal O_{\p^3}(m)\oplus \mathcal O_{\p^3}(n)\oplus \mathcal O_{\p^3}(p)\oplus \mathcal O_{\p^3}^{r-3}$$ is globally generated, hence it can be written as a quotient of $\mathcal O_{\p^3}\otimes W$ for some large vector space $W$. We obtain an immersion $i: \p^3\to {\bf G}(r, W)$ which has the property that $i^{\star} c_1(\q)=aH$, $i^{\star} c_2(\q)=bH$ and $i^{\star} c_3(\q)=cH$. Here $a, b, c$ are the elementary symmetric functions in $m, n, p$. 

Pulling back $\eqref {1m}$ from $\overline M_{0,1}(\bg, da)$ under $i$, we obtain that on $\overline M_{0,1}(\p^3, d)$:
\begin {eqnarray*}
(b\gamma_1^{(da)}&+&a^2 \gamma_7^{(da)}) ev_{1}^{\star} H^2+ (ab \gamma_2^{(da)}+ a^3\gamma_8^{(da)}) ev_1^{\star}H\cdot \kappa(H^2) + \\+(b^2\gamma_3^{(da)}&+&a^2b\gamma_4^{(da)}+a^4\gamma_9^{(da)})\cdot \kappa(H^2)^2 +(ab\gamma_5^{(da)}+c\gamma_6^{(da)}+a^3\gamma_{10}^{(da)})\cdot \kappa (H^3)\end {eqnarray*} is sum of boundaries $\Delta_{u,v}$ and boundary classes of maps with two nodes.

The independence result proved in proposition $\ref {codimension}$. {\bf B} shows that the above expression should be a multiple of $\eqref {marked}$. For example, this implies $$-\frac{1}{d}(b\gamma_1^{(da)}+a^2\gamma_7^{(da)})=ab\gamma_2^{(da)}+a^3\gamma_8^{(da)}$$ Regarding this identity as a polynomial in $a, b, c$, we immediately conclude that $$\gamma_2^{(d)}=-\frac{1}{d}\gamma^{(d)}_1, \; \gamma_{8}^{(d)}=-\frac{1}{d}\gamma_7^{(d)}.$$ The other constants can also be determined by this method: $$\gamma_3^{(d)}=\gamma_6^{(d)}=0, \gamma_5^{(d)}=-\frac{1}{d}\gamma_1^{(d)}, \gamma_4^{(d)}=\frac{1}{d^2}\gamma_1^{(d)}, \gamma_9^{(d)}=\frac{1}{d^2}\gamma_7^{(d)}, \gamma_{10}^{(d)}=-\frac{1}{d}\gamma_7^{(d)}$$ Subtracting multiple of $\eqref {marked}$ we may assume $\gamma_7^{(d)}=\gamma_8^{(d)}=\gamma_9^{(d)}=\gamma_{10}^{(d)}=0$. Moreover, the same argument will show that the only boundary terms which can appear in $\eqref {1m}$ are the classes in $(E.2)$. Therefore \begin {equation}\label{1mb} ev_1^{\star} c_2-\frac{1}{d} ev_1^{\star} c_1 \cdot \kappa(c_2)+\frac{1}{d^2}\kappa(c_1^2, c_2) -\frac{1}{d} \kappa(c_1c_2)+ (E.2) = 0\end {equation} Incidentally, this equation is enough to establish the analogous statement of proposition $\ref {codimensiona}$.{\bf B} in the case of Grassmannians. The classes in $(E)$ are generators with only two possible relations $\eqref {marked}$ and $\eqref{1mb}$. Again the description of the codimension $2$ classes on the open stratum is degree independent. 

We will now move our computation on $\overline M_{0,2}(\bg, d)$. We write equation $\eqref {1mb}$ for each of the two markings and take the difference. Using equation $\eqref {diff}$ we obtain that the following is true: 
$$ev_1^{\star}c_2-ev_2^{\star}c_2-\psi_2 \cdot \kappa(c_2)=\sum \delta'_{a,b}\Delta_{a,b}(\one, \two|c_2, \b) +\sum \delta''_{a,b} \Delta_{a,b}(\one, \two | \b, c_2).$$ We claim $\delta'_{a,b}=\delta''_{a,b}=-1$. This follows using the same technique as before, by restricting to $\overline M_{0,2}(\pr, d)\hookrightarrow \overline M_{0,2}(\bg, d)$. This time, we need to use equation $\eqref{2m}$ and the fact that the classes $\Delta_{a,b}(\one, \two|H^2, \b)$ and $\Delta_{a,b}(\one, \two |\b, H^2)$ are independent as proved in proposition $\ref {codimensiona}$. {\bf C}.

We observe that the last equation we obtained is in fact equation $\eqref {re2}$ for $\alpha = c_2(\q)$. The proof of proposition $\ref {eval}$ is now completed, in the light of remark $\ref {red}$.  

\section {The reconstruction theorem.} In this final section, we prove the reconstruction theorem $\ref{reconstructiona}$. We let $\bg$ denote the Grassmannian of dimension $2$ subspaces of a vector space $V$. We let $\langle \alpha_1, \ldots, \alpha_n\rangle_{n,d}$ denote any $n$ point degree $d$ Gromov-Witten invariant on $\bg$. The $2$ point invariants are easy to compute, they occur only in degree $1$. Henceforth, we will assume $n\geq 3$. Moreover, we assume all invariants with lower degree or fewer marked points have already been computed. We also assume that the invariants $\langle \alpha_1, \alpha_2, c_2, \ldots, c_2 \rangle_{n, d}$ are known. 

A preliminary observation is that intersections of evaluation and $\kappa$ classes are obtained from the ordinary invariants:\begin {equation}\label {kappa}\int_{\gscheme} \prod_{i=1}^{n}ev_i^{\star}\beta_i\cdot \prod_{i=1}^{m} \kappa(\alpha_i) = \int_{\overline M_{0,n+m}(\bg, d)}\prod_{i=1}^{n}ev_i^{\star} \beta_i \cdot \prod_{i=1}^{m} ev_{n+i}^{\star}\alpha_{i}.\end {equation}

Next, we claim that: \begin {equation}\label {strip}\langle \alpha_1, \alpha_2,\ldots, c_2\alpha_n \rangle_{n,d} -\langle c_2\alpha_1 , \alpha_2 \ldots, \alpha_n\rangle_{n,d}\end {equation} is sum of invariants $\langle\beta_1, \ldots, \beta_{n-1}, c_2\rangle_{n,d}$ and of invariants with fewer marked points or of lower degree. 

To explain $\eqref {strip}$, we use equation $\eqref{re2}$ for $i=n$ and $j=1$. We obtain that $\eqref{strip}$ equals $$ - \sum \int_{\gscheme} \Delta_{a,b}(S, T|c_2, \b) ev_{1}^{\star}\alpha_1\ldots ev_n^{\star}\alpha_n + \int_{\gscheme} \psi_1 \kappa(c_2) \cdot ev_{1}^{\star} \alpha_1 \cdot \ldots \cdot ev_n^{\star}{\alpha_n}.$$ Here $S$ and $T$ form a partition of $\{1, \ldots, n\}$ with $n\in S$ and $1\in T$. Let us consider the first terms in the above expression. For dimension reasons we need $a\neq 0$. Also for dimension reasons, if $b=0$, then $T$ should only have two elements, say $T=\{1,i\}$. The corresponding integral can be rewritten using equation $\eqref {kappa}$ as $\langle \alpha_1\alpha_i, \alpha_2, \ldots, \widehat \alpha_i, \ldots, \alpha_n, c_2\rangle_{n,d}$. Here the hat indicates the omission of the class $\alpha_i$. Finally, the terms with $b\neq 0$ can be expressed as usual as product of known invariants with lower degrees or fewer markings. This is a consequence of equation $\eqref{kappa}$ and the suitable analogue of the splitting axiom in Gromov-Witten theory which incorporates the $\kappa$ classes. 

For the second integral, we make use of the following folklore result (which can be derived for example from $\eqref{psisum}$): $$\psi_1=\sum \Delta_{a,b}(\{1\}, \{n-1, n\})$$ To complete the argument of $\eqref {strip}$, we use the last equation to rewrite the second integral above as a sum of known invariants with fewer markings or of lower degrees or with the class $c_2$ on the last place. Special care must be taken for the terms involving $\Delta_{a,b}(\{1\}, \{n-1, n\})$ with one component of degree $0$ carrying only two marked points. Again, we rewrite them as invariants of degree $d$ and with $n$ markings but with evaluation at the class $c_2$ on the last place. 

Equation $\eqref {strip}$ will allow us to move all powers of $c_2$ to the first position, so we can assume $\alpha_n$ is a power of $c_1$. We can then apply the Kontsevich-Manin reconstruction (\cite {KM}, 3.2.3) to finish the computation. It remains to deal with the residual terms $\langle \beta_1, \ldots, \beta_{n-1}, c_2\rangle$. In fact we can apply equation $\eqref {re2}$ a few more times to set more $\beta$'s to $c_2$. Each time, we will have to express the corresponding $\psi$ classes as sums of boundary classes. That will make use of three markings; the evaluation classes at the remaining $n-3$ markings will not be changed. Repeating this argument inductively, we can assume that all $\beta's$ except the first two are equal to $c_2$. 
These special invariants $\langle \beta_1, \beta_2, c_2, \ldots, c_2\rangle_{n,d}$ were assumed to be known. 

Finally, a dimension count, using that $\beta_1$ and $\beta_2$ can be at most $2(N-2)$ dimensional, shows that all invariants $\langle\beta_1, \beta_2, c_2, \ldots, c_2\rangle$ vanish when $\dim V>\frac{n-3}{d-2}$, with the only possible exception $d=1$. Now, the degree $1$ invariants can be easily computed from the following observation. We have a birational map obtained from the $n$ forgetful morphisms which remember only one marking: $$\pi: \overline M_{0,n}(\bg, 1)\to \overline M_{0,1}(\bg,1)\times_{\overline M_{0,0}(\bg,1)}\ldots \times_{\overline M_{0,0}(\bg,1)} \overline M_{0,1}(\bg,1).$$ This morphism is a composition of blowups over the loci where the markings coincide. It is now clear that the $n$ point degree $1$ Gromov-Witten invariants can be computed from the intersection theory of the above fibered product. The moduli spaces $\overline M_{0,0}(\bg, 1)$ and $\overline M_{0,1}(\bg, 1)$ are the flag varieties ${\bf F}(1, 3, V)$ and ${\bf F}(1, 2, 3, V)$ and the various evaluation classes on these spaces can be easily intersected. 

\begin {remark} A similar argument can be used for general varieties $X$, for primary Gromov Witten invariants with insertions in the ring generated by the Chern classes $c_1$ and $c_2$ of various vector bundles $E$ on $X$. The essential step is to prove that for all codimension $2$ classes $\alpha=c_2(E)$, the following relation holds on $\overline M_{0,n}(X, \beta)$: \begin {equation}\label {re3} \left(ev_i^{\star}\alpha-ev_j^{\star}\alpha - \psi_j\kappa(\alpha)\right)\cap [\overline M_{0,n}(X, \beta)]^{vir}=- \iota_{\star}\left(\pi_1^{\star} \kappa(\alpha)\cap \left[\Delta(\{i\}, \{j\})\right]^{vir}\right).\end {equation} Here $\iota$ denotes the inclusion of the boundary strata $\iota:\Delta(\{i\}, \{j\})\hookrightarrow \overline M_{0,n}(X, \beta)$ and $\pi_1^{\star}\kappa(\alpha)$ is the class of nodal maps on $\Delta(\{i\}, \{j\})$, the component containing the marking $i$ "intersecting" a subvariety in the class $\alpha$. 

Twisting $E$ sufficiently many times, we may assume it is globally generated, hence it induces a morphism $f:X\to \bg$ such that $f^{\star} \q=E$. Then, we pull back relation $\eqref {re2}$. We use the fact that the $\kappa$ classes pullback naturally: $f^{\star}\kappa(\alpha)=\kappa(f^{\star} \alpha),$ as the arguments of \cite {M2} show (see the reasoning leading to equation VI.3.41). An argument is required to understand how the $\psi$ classes for $X$ and $\bg$ are related under pullback. Naively, the correction terms $\psi_i^{X}-f^{\star} \psi_i^{\bg}$ are boundary divisors of nodal maps to $X$ with one component which is contracted by composition with $f$ and which contains only the marking $i$. The precise details are similar to lemma 6.6.1 in \cite {M2} where the above statement is explained for constant morphisms. The computation is slightly non-trivial, as one has to account for the virtual fundamental classes and also make use of the stack of prestable morphisms. 

Finally, equation $\eqref {re3}$ leads to a statement identical to corollary $1$. 
\addtocounter{corollary}{-1}
\begin {corollary} (bis)
Let $X$ be any smooth projective variety whose cohomology is generated by the Chern classes $c_1$ and $c_2$ of vector bundles on $X$. Fix coordinates on the small phase space and write $\bf X$ for the coordinates which are not coupled to idecomposable classes in $H^4(X)$. Then the Gromov-Witten potential $\Phi$ can be reconstructed from the inital conditions: 
$$\Phi|_{\bf X=0},\; \partial_i\Phi|_{\bf X=0},\;\text { and }\partial_{ij} \Phi|_{\bf X=0}.$$
\end {corollary}
\end {remark}

\begin {remark}
We hope our method will extend to prove more general reconstruction results. In fact, we conjecture that proposition $\ref {evala}$ holds for arbitrary classes $\alpha$ on $\bg$. We could then derive an analogue of $\eqref {re3}$ for all $X$ and all classes $\alpha$ in the ring generated by the Chern classes of holomorphic bundles (that is, via the isomorphism provided by the Chern character, for all algebraic cohomology classes on $X$). When the cohomology of $X$ is generated by the small codimension classes, we could hope for reconstruction results of the type discussed in this paper. Of course, it would be interesting to prove better reconstruction theorems which would determine, for example, {\it all} Gromov-Witten invariants of flag varieties.\end {remark}

\end {document}